\documentclass[12pt,letterpaper]{amsart}
\usepackage{palatino, euler, epic,eepic, amssymb, xypic, floatflt, microtype}
\input xy
\xyoption{all}
 \newlength{\baseunit}               % the basic unit length
 \newcount{\numlines}                % depth of picture (in number of lines)
\usepackage{amsmath, amssymb, amscd, verbatim, xspace,amsthm}
\usepackage{latexsym, epsfig, color}

\newcommand{\A}{\ensuremath{{\mathbb{A}}}}
\newcommand{\C}{\ensuremath{{\mathbb{C}}}}
\newcommand{\Z}{\ensuremath{{\mathbb{Z}}}\xspace}
\renewcommand{\P}{\ensuremath{{\mathbb{P}}}}
\newcommand{\Q}{\ensuremath{{\mathbb{Q}}}}
\newcommand{\R}{\ensuremath{{\mathbb{R}}}}
\newcommand{\F}{\ensuremath{{\mathbb{F}}}}

\newcommand{\ra}{\rightarrow}

\newcommand\Sym{\operatorname{Sym}}

\newcommand\sub{\subset}

\newcommand\Spec{\operatorname{Spec}}

\renewcommand\O{\mathcal{O}}

\newcommand\bq{\begin{equation}}
\newcommand\eq{\end{equation}}

\newtheorem{proposition}{Proposition}[subsection]
\newtheorem{theorem}[proposition]{Theorem}

\theoremstyle{remark}

\usepackage{url}

\numberwithin{equation}{section}

\newcommand{\point}{\vspace{2mm}\par \noindent
  \refstepcounter{equation}{\theequation.} } % all used to be 3mm
\newcommand{\tpoint}[1]{\vspace{2mm}\par \noindent \refstepcounter{equation}{\theequation.}
  {\bf #1. ---} }
\newcommand{\epoint}[1]{\vspace{2mm}\par \noindent \refstepcounter{equation}{\theequation.}
  {\em #1.} }
\newcommand{\bpoint}[1]{\vspace{2mm}\par \noindent \refstepcounter{equation}{\theequation.}
  {\bf #1.} }
\newcommand{\cut}[1]{}
\newcommand\hidden[1]{}

\newcommand\commentr[1]{{\color{red} \sf [#1]}}

\newcommand\highlight[1]{#1}

\newcommand{\var}{\rm{Var}}
\newcommand\wbar{\overline{w}}
\newcommand\Kbar{\overline{K}}
\newcommand\lam{\lambda}
\newcommand\la{\lambda}
\renewcommand\l{\lambda}
\newcommand{\si}{\sigma}
\newcommand\Sinf{SP}
\newcommand{\zzeta}{\mathbf{Z}}
\newcommand{\bzeta}{\boldsymbol{\zeta}}
\renewcommand{\k}{\mathbb{K}}
\newcommand{\Kvar}{K_0(\var)}
\renewcommand{\L}{\mathbf{L}}
\newcommand{\M}{\mathbf{M}}
\newcommand{\K}{\widehat{\mathcal{M}_{\mathbf{L}}}}
\newcommand{\Conf}{\operatorname{Conf}}
\newcommand{\conc}{\cdot}

\renewcommand{\cL}{{\mathcal{L}}}
\newcommand{\proj}{\mathbb{P}}
\newcommand{\cF}{\mathcal{F}}
\newcommand{\cM}{\mathcal{M}}
\newcommand{\cP}{\mathcal{P}}
\newcommand{\cQ}{\mathcal{Q}}
\newcommand{\codim}{\operatorname{codim}}
\newcommand{\FF}{\mathbb{F}}

\newcommand{\PP}{\mathbb{P}}
\newcommand{\QQ}{\mathbb{Q}}

\renewcommand{\R}{\mathbb{R}}
\renewcommand{\Z}{\mathbb{Z}}
\renewcommand{\F}{\mathbb{F}}
\renewcommand{\C}{\mathbb{C}}
\renewcommand{\P}{\mathbb{P}}
\renewcommand{\Q}{\mathbb{Q}}

\title{Discriminants in the Grothendieck Ring}

\date{August 28, 2013.}
\author{Ravi Vakil and Melanie Matchett Wood}
\address{R. Vakil,  Dept. of Mathematics, Stanford University,
  Stanford, CA 94305 USA  \newline \indent
M. M. Wood,
Dept.\ of Mathematics,
University of Wisconsin-Madison,  480 Lincoln Drive,
Madison, WI 53705 USA, and 
American Institute of Mathematics, 360 Portage Ave,
Palo Alto, CA 94306 USA} 
\email{vakil@math.stanford.edu, mmwood@math.wisc.edu}

\begin{document}
\begin{abstract}
  We consider the ``limiting behavior'' of {\em discriminants}, by
  which we mean informally the locus in some parameter
  space of some type of object where the objects have certain
  singularities.  We focus on the space of partially labeled points on
  a variety $X$, and linear systems on $X$.  These are connected ---
  we use the first to understand the second.  We describe their
  classes in the Grothendieck ring of varieties, as the number of
  points gets large, or as the line bundle gets very positive.  They
  stabilize in an appropriate sense, and their stabilization is given
  in terms of  motivic zeta values.  Motivated by our results, we
  conjecture that the symmetric powers of geometrically irreducible
  varieties stabilize in the Grothendieck ring (in an appropriate
  sense).  Our results extend parallel results in both arithmetic and
  topology. We give a number of reasons for considering these
  questions, and propose a number of new conjectures, both arithmetic
  and topological.
\end{abstract}
\maketitle
\setcounter{tocdepth}{1} % this just includes subsections
\tableofcontents

\section{Introduction}

\label{s:int}We study the classes of 
discriminants (loci in a moduli space of objects with specified
singularities) and their complements in the Grothendieck ring of
varieties, focusing on the cases of moduli of hypersurfaces and
configuration spaces of points.  
The main contributions of this paper are two theorems (Theorems~\ref{T:Bjorn} and~\ref{T:IntroPointsLimits}) and one conjecture (``motivic stabilization of symmetric powers'',
Conjecture~\ref{c:mssp}). 

\noindent {\bf I.}
{\em (Theorem~\ref{T:Bjorn}, the limiting motive of the space of
    hypersurfaces with a given number of  singularities, \S \ref{S:Hyper})} If $\cL$ is an ample line bundle on a smooth variety $X$, we
  show that the motive of the subset of the linear system
  $|\cL^{\otimes j}|$
 consisting of divisors with
  precisely $s$ singularities
  (normalized by $|\cL^{\otimes j}|$), tends to a limit as $j \rightarrow
  \infty$ (in the completion of the localization of the
Grothendieck ring at $\L := [ \A^1 ]$), given  explicitly in terms of the motivic zeta function
  of $X$.  
  
\noindent {\bf II.}  {\em (Conjecture~\ref{c:mssp}, motivic
  stabilization of symmetric powers, \S \ref{S:MSSP})}
We conjecture that if $X$ is geometrically irreducible, then the ratio
$[\Sym^n X] /
\L^{n \dim X}$ tends to a limit.  This is an algebraic version of the
Dold-Thom theorem, and is also motivated by the Weil conjectures.  We
give a number of reasons for considering this conjecture.

\noindent {\bf III.}  {\em (Theorem~\ref{T:IntroPointsLimits}, the
  limiting motive of discriminants in configuration spaces, \S \ref{s:cs})} We show
that if $X$ is geometrically irreducible and satisfies motivic
stabilization ({\bf II}, e.g.\ if $X$ is stably
rational, see Motivation~\ref{motivation}(i)), then the motive of strata (and their closure) of
configurations of points with given ``discriminant'' (clumping of
points) tends to a limit as the number of points $n \rightarrow
\infty$, and (more important) we describe the limit in terms of
motivic zeta values.  In the case of $s$ multiple points, the result
is the same as that of {\bf I}, except the expression in terms of
motivic zeta functions is evaluated at a different value  (see
Theorem~\ref{T:IntroPointsssing} and \eqref{eq:blimp}).  The
reliance on the motivic stabilization conjecture can be removed by specializing
to Hodge structures, where the analogous conjecture holds (see
Motivation~\ref{motivation}(ii)), or by working with generating series
(Theorem~\ref{T:Kform}).
 
These results are motivated by a number of results in number theory
and topology (including, notably, stability/stabilization theorems), and they generalize
analogues of many of these statements.  (An elementary motivation is
an analogue of both {\bf I} and {\bf III} for $X= \Spec \Z$: the
probability of an integer being square free is $1 / \zeta(2)$.  One
has to first make sense of the word ``probability'' as a limit, then
show that the limit is a zeta value.  These features will be visible
in our arguments as well.)  Our results also support Denef and
Loeser's motto  \cite[l.~1-2]{Denef-Loeser2004}:  ``rational generating series occurring in arithmetic
geometry are motivic in nature''.

Our results suggest a large number of new conjectures in arithmetic,
algebraic geometry, and topology that may be tractable by other means.
We label certain of these smaller new conjectures by letters~A
through~H to draw attention to them among the many conjectural
statements we make.  The reason for stating many of these conjectures
is not necessarily an expectation that they will be true, but because
either a proof or a counterexample should provide significant new
insight.  Hence throughout this paper, {\em conjecture} should perhaps
be interpreted as {\em conjecture/speculation/question}.  Since we
made these conjectures publicly in the first draft of our paper, they
have motivated a significant amount of research by several authors
into the structure that our conjectures highlight, e.g \cite[Cor.~3]{Church12} and \cite{KM, KM2,
  Tom, Tom2}.  Some of the conjectures have been proven in special cases,
some have been proven entirely, and some have been disproven.  One
important feature of our work is that it has brought to light
structure in certain questions of homological stability, that despite
being in a very active area of research, had been previously
overlooked.

We now describe these results in more detail and context.  We first set
notation and review background about the Grothendieck ring.
Thereafter, the discussion  of {\bf I}, {\bf II}, and {\bf III} can be read largely in any order.

\bpoint{Notation and background:  The Grothendieck ring of
  varieties}
\label{s:Gring}\label{s:notdef}Throughout, \highlight{$\k$} is a field.  A {\em variety} is
a reduced separated finite type $\k$-scheme; \highlight{$X$} will
usually denote a variety, and \highlight{$d$} will be its
dimension.

The Grothendieck ring of varieties \highlight{$\cM := K_0(\var_{\k})$}
is defined as follows.  As an abelian group, it is generated by the
classes of finite type $\k$-schemes up to isomorphism.  The class of
a scheme $X$ in $\cM$ is denoted \highlight{$[X]$}, but we often drop the
 brackets for convenience.  The group relations are generated by
the following: if $Y$ is a closed subscheme of $X$, and $U$ is its
(open) complement, then $[X] = [U] + [Y]$.  In particular, taking
$Y=X^{red}$, we have $[X] = [X^{red}]$, so nilpotents play no role in
our discussions.  The product $[X][Y] :=[X\times_{\k} Y]$ makes $\cM$
into a commutative ring, with $[\Spec \k]$ as unit.  

Any morphism $\phi$ from $\cM$ to another ring is called a {\em
  motivic measure}.  Here are two important examples. (i) If $\k = \FF_q$,
there is a {\em point counting map} \highlight{$\#: \cM \ra \Z$}.

(ii) For  convenience, we call the Grothendieck
group of mixed Hodge structures {\em virtual Hodge structures}. 
This is the same as the sum (over all weights) of the
Grothendieck group of $\Q$-Hodge structures.
 If $\k = \C$, there is a map from $\cM$ to the group of virtual Hodge
 structures,  defined by taking each 
variety $X$ to $\sum_k (-1)^k [H^k_c(X, \Q)]$.   This descends to a motivic measure
\highlight{$HS$} from $\cM$ to the group of virtual Hodge structures:
 given complementary subsets $Y, U \subset X$ as in two paragraphs previous, the long exact sequence for cohomology with compact
support respects mixed Hodge structure.  This motivic measure
specializes further to the {\em Hodge-Deligne polynomial} $e: \cM
\rightarrow \Z[x,y]$, where for a variety $X$,
$$
e(X) = \sum_k (-1)^k h^{p,q}(H^k_c(X)) x^p y^q.
$$
(This has also been called the $E$-polynomial, the
virtual Hodge polynomial, the Serre polynomial, and the Hodge-Euler
polynomial.)  
If $X$ is smooth and proper, $e(X)$ determines the Hodge
numbers on each of the cohomology groups, and in particular the Betti
numbers $h^i(X)$.  
\cut{(See \cite{DeligneHodgeII, DeligneHodgeIII} for
background on mixed Hodge structures on cohomology; \cite{DanKhov86}
for the definitions using compactly supported cohomology based on
Deligne's construction of mixed Hodge structures on relative
cohomology \cite[8.3.8]{DeligneHodgeIII}; and \cite{PetersSteen08} for
proofs of the properties of the mixed Hodge structures. )}

\epoint{Principle: Occam's razor for Hodge
  structures} \label{s:occraz}We point out a well-known (if vague)
principle.  For a variety $X$ that is not smooth and proper, the
virtual Hodge structure does not determine the Hodge structures on
each $H^k_c(X)$ (and similarly the Hodge-Deligne polynomial does not
determine the Hodge numbers $h^{p,q}(H^k_c(X))$) because the
contributions from different $k$ are mixed.  But in many cases
%(especially if $X$ is smooth or proper) 
there is a simplest Hodge
structure on all the $H^k_c$'s compatible with the virtual Hodge
structure, and it is reasonable to wonder if the Hodge structure on
the $H^k_c$'s is this simplest possibility in some simple cases.
Similarly, if the virtual Hodge structures stabilize in some sense in
some family of examples, one should expect that this arises because
the actual Hodge structures on the $H^k_c$'s (and hence the ``compact
type'' Betti numbers) also stabilize.

\epoint{Inverting $\L$} We denote the ``Lefschetz motive'' $[\A^1]$ by
\highlight{$\L$}.  There are many reasons to consider the localization
\highlight{$\cM_\L$} (including motivic integration; possible rationality of the motivic zeta
function, \S \ref{q:mzfrat}; and the desirability of a homotopy axiom).  This paper suggests additional
reasons.  The motivic measures $HS$ and $\#$ clearly extend to
$\cM_\L$.

\epoint{Completion with respect to the dimensional filtration} \label{s:dimfil}The Grothendieck ring $\cM$  is filtered by the subgroups generated by varieties
of dimension at most $d$, as $d$ varies.  This {\em dimensional filtration}
clearly extends to $\cM_{\L}$.  Let \highlight{$\K$} be the completion
of $\cM_{\L}$ with respect to the dimensional filtration.  As explained
in \cite[\S 1.5]{Bourqui2010} , $\K$
inherits not just a group structure, but also a ring structure. The
ring $\K$ was originally introduced by Kontsevich \cite{Klecture95} as
the ring where values of motivic integrals lie, in his theory of
motivic integration (see \cite{DL99, DL01, Loo02}).
Note that the motivic measure $HS$ extends to this
completion, after suitably extending the codomain.  But
the point-counting motivic measure does not
 --- the point counting map $\#: \cM_{\L} \rightarrow \Q$
is not continuous (consider the sequence
$(2q)^n\L^{-n}$).

\epoint{Observation} The symmetric product $\Sym^n X$ will be central
to us.  If $X$ is not quasiprojective, then $\Sym^n X$ might not be a
variety.  However, $[\Sym^n X]$ may be interpreted as an element of $\K$, by
\cite[Thm.~1.2]{ekedahl} (note that $\Sym^n X$ is represented by an
algebraic space), which will suffice for our purposes; we will
hereafter pass over this technical point without comment.  We note for
future reference that (i) if $X$ is rational, then $[\Sym^n X]$ is
invertible in $\K$, and (ii) if $\k=\C$ then $HS(\Sym^n X)$ is
invertible in $HS(\K)$ (for any $X$).\label{o:Syminv}

\epoint{The motivic zeta function}
\label{s:mzf}Let \highlight{$\zzeta_X(t) :=\sum_{n\geq 0} [\Sym^n X]t^n
  \in \cM[[t]]$} be the {\em motivic zeta function} (defined by
Kapranov, \cite[(1.3)]{Kapranov2000}).  If $\k =\F_q$, then
$\#$ sends $\zzeta_X(t)$ to the Weil zeta function $\zeta_X(s)$, where
$t=q^{-s}$.   Motivated by this, we define (for {\em any} $\k$) \highlight{$\bzeta_X(m) :=
\zzeta_X(\L^{-m})$}.  We use bold fonts for both $\bzeta$ and $\zzeta$
in order to distinguish them from the Weil zeta function(s).

 Let \highlight{$\Sym^n_{[s]} X \subset \Sym^n X$}
(not to be confused with $\Sym^n_{s} X \subset \Sym^n X$, to be
defined in \S \ref{d:syms})
 be the
locally closed subset of $\Sym^n X$ consisting of unordered $n$-tuples
of points supported on exactly $s$ (distinct, geometric) points.  Let
\highlight{$\zzeta^{[s]}_X(t) = \sum_n [\Sym^n_{[s]} X ] t^n$}, so
$\zzeta_X(t)=\zzeta^{[0]}_X(t)+\zzeta^{[1]}_X(t)+\zzeta^{[2]}_X(t)+\dots$.
It is straightforward to write $\zzeta^{[s]}_X(t)$ in terms of $\Sym^n
X$'s and rational functions of $t$.  For example,\hidden{details are
  in discnotesmay0112.tex, just before the ``part'' on the
  $L$-conjecture}\begin{equation}\label{eq:012}\zzeta^{[0]}_X(t) = 1,
  \quad \zzeta^{[1]}_X(t) = \frac {t} {1-t} X,   \quad \text{and}\end{equation}
$$\zzeta^{[2]}_X(t) = \frac {  t^2}
{1-t^2} \Sym^2 X+ \frac { t^3} { (1-t^2) (1-t)} X^2 -  \frac {t^2} {
  (1-t^2) (1-t)} X.$$
Define \highlight{$\bzeta^{[s]}_X(m) := \zzeta^{[s]}_m(\L^{-m})$}.

\epoint{Rationality of the motivic zeta
  function} \label{q:mzfrat}Kapranov \cite[Rem.~ 1.3.5]{Kapranov2000} asked if the motivic zeta
function is rational, given that its specialization $\zeta_X(t)$
{\em is} rational.  
(This question is related to {\bf II}, see \S \ref{c:mssp}(iv).) 
Further evidence is that
the motivic measure to Hodge structures $HS(\zzeta_X(t))$  is rational, which was shown
by Cheah (actually predating the definition of the
motivic zeta function).   We note that some care is necessary to say
what is meant by rationality, see \cite[\S 2]{LL04}.

\tpoint{Theorem (Cheah, \cite{Cheah94}, see also \cite{Cheah96})}
{\em  Suppose $X$ is a complex variety. 
Then 
$$
HS(  \zzeta_X(t) ) = \prod_{i=0}^{\infty}  \left( { (1-t)^{ H^i(X)}
} \right)^{(-1)^{i+1}}
$$
where  if $V$ is a mixed Hodge structure, then $(1-t)^{[V]}$ is 
interpreted as $\sum_{j=0}^{\infty}    (-1)^j [ \wedge^j V] t^j$, and $[\cdot
]$ indicates the class in virtual Hodge structures.}\label{t:cheah}

Cheah's argument deals only with the specialization to the Hodge-Deligne
polynomials (actually an enrichment of this), but the proof can be
adapted to establish Theorem~\ref{t:cheah}.

\point \label{s:ll}However, Larson and Lunts \cite{LL03,LL04} showed that $\zzeta_X(t)$ is
{\em not} always rational in $\cM[[t]]$.  But an important (if vague) question
remains: where between $\cM$ and the special motivic measures of
point-counting or Hodge structures is the motivic zeta function
rational?  In particular, the argument of Larson and Lunts does not
apply to $\cM_{\L}$. 

\tpoint{Conjecture \cite[Conj.~ 7.5.1]{Denef-Loeser2004}} {\em \label{q:zrat}The generating series $\zzeta_X(t)$ is rational in $\cM_{\L}[[t]]$.}

\bpoint{Moduli of Hypersurfaces}

Our main result on hypersurfaces is the following, proven in Section~\ref{S:Hyper}.

\tpoint{Theorem}\label{T:Bjorn}
{\em Let $X$ be a smooth projective variety of pure dimension $d>0$
with an ample line bundle $\cL$. 
Let $H^0(X,\cL^{\otimes j})^s$ be the constructible subset of
$H^0(X,\cL^{\otimes j})$ corresponding to divisors on $X$ with exactly
$s$ singular geometric  points.
Then }
\begin{equation}\label{eq:Bjorn}
\lim_{j \ra \infty} \frac{  [H^0(X,\cL^{\otimes j})^s]}{[
  H^0(X,\cL^{\otimes j})]}= \frac  {\bzeta^{[s]}_X(d+1)   }
{\bzeta_X(d+1)  } \quad \quad \text{(in $\K$).}
\end{equation}

For example (using \eqref{eq:012}),   the motivic (limiting) probability of a divisor
being smooth (i.e.\ $s=0$) is $$1/\bzeta_X( d+1),$$ and  the
motivic
(limiting) probability of a divisor having precisely one singularity
($s=1$) is  
$$
\frac{ X  \L^{-(d+1)}  }   {1- \L^{-(d+1)}}  \cdot \frac 1
{\bzeta_X( d+1) }.
$$

\epoint{Remarks} \label{r:Bjorn} \quad

(i) Note that the limiting motivic
density (the right side of \eqref{eq:Bjorn}) is independent of $\cL$. 

(ii) Even to establish the result for
$\k=\C$, the argument requires the use of finite fields (see
Lemma~\ref{L:SymInd}).  

(iii) The hypothesis of projectivity can be
weakened to quasiprojectivity by taking appropriate care in defining
the motivic probability (cf.\ \S \ref{m1}).  We leave this variation to the interested
reader (see Motivation~\ref{m1}).

(iv) If $s>1$, $H^0(X,\cL^{\otimes j})^s$ will in general not be
locally closed.  For example, consider $s=2$ and $d=2$, in the
neighborhood of a curve with a tacnode and a node. 

(v)  If instead of $s$
singular points, we require $s$ singular geometric points that are
$A_1$-singularities, then the corresponding locus is
locally closed (not just constructible), and the limit again exists,
and can be interpreted as follows.  Let ${\mathbf{P}}_d$ be the motivic
probability that a singular point $p$ of a divisor on a smooth
$d$-fold  $X$ passing through $p$ is an $A_1$-singularity.  (This can be made precise in the
obvious way, by considering $\Sym^2( \Omega_X|_p)$.)
Then multiply the right side of \eqref{eq:Bjorn} by $({\mathbf{P}}_d)^s$.  We omit
the justification.  However, see Conjecture \ref{c:b} below.\label{r:} 

(vi)
A simple extension of the argument yields the following result: the
motivic (limiting) probability that a section of $\cL^{\otimes j}$ has
no $m$-multiple points is $1/ \bzeta_X\left( \binom{d+m-1}{d} \right)$.  (An $m$-multiple point is a singular point of
multiplicity at least $m$ --- the defining equation vanishes to order
at least $m$.)

(vii) A variation of the proof of Theorem~\ref{T:Bjorn} 
(see Remark~\ref{r:sordered})
yields the
following.  Let $H^0(X,
\cL^{\otimes j})^{\text{$s$ ordered}}$ be the space of sections of
$H^0(X, \cL^{\otimes j})$ along with a choice of $s$ {\em ordered}
(disjoint) singular
(geometric) points.  Then 
$$\lim_{j \ra \infty} \frac{  [H^0(X,\cL^{\otimes j})^\text{$s$ ordered}]}{[
  H^0(X,\cL^{\otimes j})]}= 
\frac {  [  X^s \setminus \Delta]}  { \bzeta_X(  {d+1})}  \left(
  \frac { 1 / {\L^{d+1}}  } { 1 - 1 / {\L^{d+1}}   } \right)^s
 \quad \quad \text{(in $\K$)}
$$
where $\Delta \subset X^s$  is the ``big diagonal''.

We give three motivations for Theorem~\ref{T:Bjorn}.

\epoint{First Motivation:  Poonen's probability for a hypersurface to be
  smooth}
\label{m1}Poonen's ``Bertini Theorem over finite fields'' \cite[Thm.~1.1]{Poonen04} is (informally) the following.  Suppose $X
\subset \proj^N$
is a smooth projective variety over $\F_q$ of dimension $d$.
(Poonen  states his result more generally in the quasiprojective
case, and ours can be so extended as  well, see Remark~\ref{r:Bjorn}(iii).)  As the base field is finite, one can make sense of the
probability $p_j$ of a hypersurface of degree $j$  (defined over
$\F_q$) intersecting $X$
along a smooth (codimension $1$) subvariety.  Poonen shows that
$\lim_{j \rightarrow \infty} p_j$ exists, and equals $1/\zeta_X(d+1)$,
where $\zeta_X$ is the Weil zeta function.

Theorem~\ref{T:Bjorn} in the case $s=0$ is the motivic analogue of
Poonen's result.  However, this case of Theorem~\ref{T:Bjorn} neither
implies nor is implied by \cite[Thm.~1.1]{Poonen04}; the limits
taken in both cases are not compatible, because the dimensional
filtration has no relation to point-counting (\S \ref{s:dimfil}).
Furthermore, the methods of proof are unrelated (except at a very
superficial level).  Based on Theorem~\ref{T:Bjorn}, it is reasonable
to conjecture the following.

\tpoint{Conjecture~A}
{\em  Suppose $X\subset \proj^N$
is a smooth, dimension $d$,
variety over $\F_q$. 
Let $p_j^{[s]}$ be the 
probability of a hypersurface of degree $j$ (defined over $\F_q$) intersecting $X$
along a subvariety with exactly $s$ 
non-smooth geometric points.  Then 
$$\lim_{j \rightarrow \infty} p_j^{[s]} = 
\frac  {\zeta^{[s]}_X(d+1) }  {\zeta_X(d+1)},$$
where $\zeta^{[s]}_X$ is defined analogously to $\bzeta^{[s]}_X$ (\S \ref{s:mzf}).
}

The case of $\proj^1$ follows from
Theorem~\ref{T:IntroPointsssing}(b), and the heuristic given after the
statement of Theorem~\ref{T:IntroPointsssing} suggests this result for
smooth curves ($d=1$) in general.

\epoint{Second Motivation:  motives of Severi varieties}
Severi varieties are closed subsets of linear systems on a smooth
projective surface $(X, \cL)$ with a fixed number of
singularities. G\"ottsche's conjecture (now a theorem, see \cite{T12}  and \cite{kst}) states that the
degrees of these varieties have strong structure: if the line bundle
is sufficiently ample, then the degree can be read off from a
universal formula involving four universal generating series (two of
them quasimodular forms), and the four numerical invariants of $(X,
\cL)$.  (See \cite{kst} for a more precise statement.)

Theorem~\ref{T:Bjorn} states that not only does the degree of the
Severi variety have a strong structure (related to modular forms),
the motive does as well (related to zeta functions).  (If one wishes
to restrict to nodal curves, i.e.\ $A_1$-singularities, as is usually the case for Severi
varieties, see Remark~\ref{r:}(v).  We conjecture similar
structure if one considers more general singularity types, in vague
analogy with \cite{Li-Tzeng2012}, extending the proof \cite{T12} of
G\"ottsche's conjecture, or \cite{r}, extending the proof \cite{kst}.)

\epoint{Third Motivation: Vassiliev's work on topology of
  discriminants and their complements} \label{s:v}Vassiliev's fundamental work on
topology of discriminants or discriminant complements (in a broad
sense, encompassing different singularity types) can be seen as
philosophical and direct
 motivation for this work.  We left it for last as it
leads to further questions and conjectures.  We briefly summarize some
of Vassiliev's work; see \cite{V95} and the references (by Vassiliev)
therein for more.  Vassiliev determines the 
cohomology ring of the space of holomorphic functions in
$d$ complex variables giving a smooth divisor, 
 \cite[Thm.~1]{V90}.  His argument works without change to
apply to {\em algebraic} (i.e.\ polynomial) functions in $d$ complex
variables.  This predicts the same Betti numbers as the application of
Occam's Razor~\ref{s:occraz} to Theorem~\ref{T:Bjorn} in the special
case $X = \C^d$, $\cL = \mathcal{O}$ (see Remark~\ref{r:Bjorn}(iii)),
which suggests that if $S$ is the space of algebraic functions with
smooth divisor, $h_i(S)=0$ for $i \neq 0, 1$, and $h_0(S) =
h_1(S) =1$.

Thus motivated, we conjecture the following.

\tpoint{Conjecture~B} {\em Suppose $X$ is a smooth complex 
  variety, and let $\cL$ be an ample line bundle on $X$.  Let $X^s_j$
  be the space of sections of $\cL^{\otimes j}$ that vanish on a
  divisor of $X$ singular at precisely $s$ points, each of
  which is an $A_1$-singularity (so $X^0_j$ is the space of sections
  whose vanishing scheme is a smooth divisor of $X$).  Let $Y^1_j$ be the space
  of sections that vanish on a divisor of $X$ singular at one point.
  Then the rational homology type of $Y^1_j$ and each $X^s_j$
  stabilizes (i.e. for
 every $i$, we have that the singular homology groups $h_i(Y^1_j)$ and $h_i(X^s_j)$ stabilize for
 $j\gg_i 0$), and each limit is independent of
  $\cL$.\label{c:b}}

See Remarks~\ref{r:Bjorn}(iv) and (v)  for motivation.  The case most of
interest is $s=0$, and even when $X$ is a projective or
affine variety this case is not clear. 
(One might hope that Vassiliev's arguments can be extended to this case.
Note in particular that his constructions can be algebraized, and that
the spectral sequence used in his proof degenerates at $E^1$, see 
\cite[p. 212]{V95}.)

Also, even in the case of $X = \C^d$, this conjecture is not clear
(except for $s=0$, which is Vassiliev's result).   Since the authors
publicly made Conjecture B in the first draft of this paper, Tommasi
\cite{Tom2} has proven that Conjecture B for $X^s_j$ is true in the
cases where $X=\P^n$ and $\L=\O(1)$ and $s\leq 2$, and moreover has
found the stable homology in these cases.  One can also give an
analogous form of Conjecture~B for divisors of $X$ with no singularity
of multiplicity $m$ (or exactly one such).

It would be interesting not just to know that the limits in Conjecture~B (and its variant with $m$-fold points) exist, but to actually describe the limit rational homology type (i.e., the Poincare series), and in particular,
compare these limits to the limit motives given in Theorem~\ref{T:Bjorn}.
 In
particular, one could hope there is some imprecise dictionary
between motivic zeta-values of $X$ and rational homology types built in some
way out of $X$.  
Somewhat more precisely, when one sees a motivic
zeta-value in a limiting formula for some geometric problem, one
might expect to see the corresponding rational homology type in the
stabilization.

\tpoint{Vague question} {\em What rational homology type does $1/\bzeta_X(N)$
  correspond to?}\label{vg}

This vague question suggests explicit questions. For example, in
\cite[p.~1066 (2)]{T96}, Totaro gives  two  simple
complex projective manifolds, $X = \proj^1 \times \proj^2$ and $Y =
\proj_{\proj^2} (\mathcal{O} \oplus \mathcal{O}(1))$, with $[X] = [Y]$ (and so $X$ and
$Y$ have the same Betti numbers), such that the space of $3$ distinct
points on $X$ has different Betti numbers than the space of $3$
distinct points on $Y$.  In light of Theorem~\ref{T:IntroPointsLimits}(a) below (with
$\lambda = \varnothing$), Vague Question~\ref{vg} leads one to ask if this
pathology goes away as the number of points gets large.  More
precisely, does the cohomology of the space of $n$ distinct points on
$X$ agree with that of $n$ distinct points on $Y$ in any particular
degree as $n \rightarrow \infty$?    (Totaro pointed out to us that
\cite{ft} has a different description of the Betti numbers of configuration
spaces that may be better to approach this problem.)

We also conjecture that  when the Hodge-Deligne series of the limiting motive is finite, then the limit of the corresponding
Poincare series is finite.
As an important explicit 
  example, we have the following.

\tpoint{Conjecture~C} {\em For a smooth, complex variety $X$, if $e(\bzeta_X(d+1))$ is a polynomial in $x^{-1}, y^{-1}$ (e.g.\ if $X$ is 
projective with no odd degree cohomology), then for $i$ sufficiently large, and $j$ sufficient large depending on $i$,
we have $h_i(X^0_j)=0$.
}

%One can also make versions of Conjecture~C for $X^s_j$, $Y^1_j$, or the no or $1$ $m$-multiple singularities version.

Under the hypotheses of Conjecture~C, it is natural to wonder whether an application of Occam's
Razor~\ref{s:occraz} to the limit motive gives the correct limit Betti numbers.
For $X=\A^d$, the work 
 of
Vassiliev mentioned in \S \ref{s:v} implies that the answer is yes.
For  $\P^1$, 
\cite[Prop.~4.5]{Church12v1} can be used to show that the answer
is no for $X^0_j$, but yes for $X^0_j/\C^*$.  (We thank T.~Church for
explaining this to us, \cite{churchpc}.)\cut{see his
  email Jul 23 2012}
%and for $X=\A^d$ a positive answer is further suggested by Vassiliev's
%result above.  
As just one open example for $X^0_j/\C^*$, we highlight the case of plane curves.

 \tpoint{Conjecture~D} {\em Let $X^0_j / \C^*$ be the space of 
degree $j$ 
%homogeneous polynomials in $\C[x_0,x_1,x_2]$ that give 
smooth
  projective plane curves (a  quasiprojective
   manifold of complex dimension $\binom {j+2} 2-1$).
Then $\lim_{j\ra\infty} h_i(X^0_j)=1$  for $i=0,3,5,10$
%if $i=0,1,3,8,10,11$ FOR SPACE OF SECTIONS
and
$\lim_{j\ra\infty} h_i(X^0_j)=0$ otherwise. }

Since the authors publicly made this conjecture in the first draft of
this paper, Tommasi \cite{Tom} has proven that Conjecture F is
incorrect, and in fact $\lim_{j\ra\infty} h_i(X^0_j)=1$ for
$i=0,3,5,8$ and $\lim_{j\ra\infty} h_i(X^0_j)=0$ otherwise.  The
stable degree $8$ homology class has weight $10$.  Furthermore, Tommasi
\cite{Tom} has found the stable homology when $\P^2$ is replaced by
any $\P^n$.  Though there has previously been work to compute homology
groups of spaces of smooth hypersurfaces in small cases, the key
insight from our paper is that there is nicer structure if one studies
the stable homology groups.

We remark that the fundamental group of $X^0_j / \C^*$ (parametrizing
smooth plane curves) was computed by
L\"onne, \cite[Main Thm.]{lonne}, and also behaves well as $j
\rightarrow \infty$.\cut{Zariski's earlier AJM paper and Annals paper
  may be relevant too --- but I don't know what Poincare group means
  to Zariski.  -R}

\bpoint{Motivic stabilization of symmetric powers} \label{s:mssp1}

\tpoint{Conjecture (Motivic stabilization of symmetric powers)}
{\em Suppose $X$ is a geometrically irreducible variety of dimension $d$.
Then 
the limit
$
\lim_{n\ra\infty} {\left[ \Sym^n X \right]}/{\L^{dn}}
$
exists in $\K$.}\label{c:mssp}

(To see the necessity of the geometric irreducibility hypothesis, consider
the case where $X$ is two points, or see Motivation~\ref{motivation} below.)  If
Conjecture~\ref{c:mssp} holds for $X$, we say that \highlight{MSSP (or
  motivic stabilization of symmetric powers) holds for $X$}.  If
$\phi: \K \ra \phi(\K)$ is a continuous ring homomorphism (extending a
continuous motivic measure $\phi: \cM_{\cL} \rightarrow \phi(\K)$), we say
\highlight{MSSP${}_\phi$ holds for $X$} if
$$
\Sinf_\phi(X) := \lim_{n\ra\infty} \frac{\left[ \Sym^n X \right]}{\L^{dn}} \textrm{ exists in
  $\phi(\K)$,}
$$
where we abuse notation by using $[Z]$ to denote $\phi([Z])$.
Our use of the topology notation  \highlight{$\Sinf$(X)}  for infinite
symmetric product  is motivated
by the Dold-Thom theorem, see Motivation~\ref{motivation}(v) below.
(It may be suggestive to write 
$[ \Sym^n X] / [\L^{dn}]$ as $[\Sym^n X] / [\Sym^n \A^{\dim X}]$,
using Proposition~\ref{p:A}.)

\epoint{Motivation}  \label{motivation}We give a number of motivations for considering Conjecture~\ref{c:mssp}.

(i) We show that Conjecture~\ref{c:mssp} holds (or fails) on stable
birational equivalence classes (Proposition~\ref{p:A} combined with
Proposition~\ref{p:B}), and for curves with a rational point
(Proposition~\ref{p:C}).  In particular, as Conjecture~\ref{c:mssp}
clearly holds for a point, it holds for all stably rational varieties.

(ii) Conjecture~\ref{c:mssp} is true upon specialization to Hodge structures (i.e.,
MSSP${}_{\text{HS}}$ holds, for all $X$), which can be shown from Theorem~\ref{t:cheah}. 

(iii)  The analogue for point-counting also holds.  (Unlike (ii), this is
only an analogy:
point-counting is not compatible with the completion with respect to
the dimensional filtration, see \S
\ref{s:dimfil}.)  More precisely, if $X$ is a geometrically
irreducible variety over a finite field, then $\lim_{n \rightarrow
  \infty} \frac {\# \Sym^n X} { q^{dn}}$ exists.  This is
because, by the Weil conjectures, the generating function for the Weil
zeta function $\zeta_X(t)$ has as its denominator a polynomial whose
smallest root is $1/q^{d}$ (corresponding to the fundamental
class of $X$), and this root appears with multiplicity $1$. 

(iv) Related to Conjecture~\ref{q:zrat} on whether the motivic
zeta function $\zzeta_X(t)$ is rational (upon localization by $\L$), one
is led to ask whether, in a suitable sense, the denominator of $\zzeta_X(t)$
has a unique smallest root (in the sense of dimension), $\L^{- d}$ (corresponding to the ``fundamental class of $X$'', in further
analogy to the Weil conjectures).  Suitably interpreted, this would imply
Conjecture~\ref{c:mssp}.  

(v) A topological motivation is the Dold-Thom theorem \cite{dt}\cut{ (see for
example [Hatcher2002, \S 4.K])}, and more basically that the homotopy
type of $\Sym^n X$ has a limit $\Sinf (X)$, where $X$ is a topological
space (see, for example,  \cite[\S 2]{ccmm} for more discussion).  
 If $\k=\C$, then  Dold-Thom
implies that $h_i(\Sym^n X, \Q)$ stabilizes as $n\ra \infty$. 
If further $X$ is smooth, 
then Poincar\'e duality holds for $\Sym^n X$ (with $\Q$-coefficients),
because $\Sym^n X$ is the coarse moduli space for the orbifold (smooth
Deligne-Mumford stack) $X^n / \mathfrak{S}_n$.
The quotient by $\L^{dn}$ in the
statement of Conjecture~\ref{c:mssp} arises because then, by Poincar\'e duality,
$h^{2nd-i}_c(\Sym^n X)$ stabilizes as $n\ra \infty$, and the weight $-i$ piece of $e([\Sym^n X]/L^{dn})$
is the weight $2nd-i$ piece of $h^{*}_c(\Sym^n X)$.\cut{yyy I'll tweak
  this further after hearing responses to my emails.}

(vi)  Kimura and Vistoli have given analogous conjectures for Chow
groups, notably their Weak Stabilization Conjecture
\cite[Conj.~2.6]{kv} (true for curves, \cite[Prop.~2.9(a)]{kv}) and their Strong
Stabilization Conjecture \cite[Conj.~2.13]{kv} (true for pointed
curves of genus up to $4$, \cite[Cor.~2.19]{kv}, and with motivation for all curves, \cite[Rk.~2.20]{kv}).

(vii)  The statement of Conjecture~\ref{c:mssp}
contradicts each of two  well-known questions (or conjectures),
Conjectures~\ref{c:cut} and~\ref{c:L} below, as shown by D. Litt,
\cite{Litt12}.

\tpoint{Piecewise Isomorphism Conjecture (Larsen-Lunts,
  \cite[Qu.~1.2]{LL03}; see also 
\cite[Assertion 1]{LiuSebag10})}
{\em \label{c:cut}If $X$ and $Y$ are varieties with $[X] = [Y]$ in $K_0(\var_{\k})$, then 
we can write $X = \coprod_{i=1}^n X_i$ and
$Y = \coprod_{i=1}^n Y_i$ with $X_i$ and $Y_i$ locally closed, and
$X_i \cong Y_i$ ($X$
and $Y$ are ``piecewise isomorphic'').}

Liu and Sebag have proved Conjecture~\ref{c:cut} if $\k$ is
algebraically closed of characteristic $0$, when $X$ is a smooth
connected projective surface \cite[Thm.~4]{LiuSebag10}, or when $X$
contains only finitely many rational curves \cite[Thm.~5]{LiuSebag10}.

\tpoint{Conjecture (well-known, see for example \cite[\S 3.3]{Denef-Loeser2004})}
{\em The element $\L$ is not a zero-divisor.  Equivalently,  the
localization $\cM
\rightarrow \cM_{\L}$ is an injection.\label{c:L}} 

This
is more a question than a conjecture.  The real (if vague) question is:
``what information, if
any, is lost by localizing by $\L$?''  (It is known that  $\cM$ is not
an integral domain, 
see \cite[Thm.~1]{Poonen02}, \cite[Ex.~6]{Kollar05}, 
\cite[Thm.~22]{Naumann07}.)
In light of Larsen and Lunts' counterexample to the rationality in
general of
the motivic zeta  function, \S \ref{s:ll}, Conjecture~\ref{c:L}
contradicts Conjecture~\ref{q:zrat}.

\bpoint{Configurations of points on varieties}
\label{s:conspa}

For a partition $\lambda$ of $n$, let $w_\lambda( X)$
be the locally closed subset of $\Sym^n X$ that is the locus of points
which occur with multiplicities precisely $\lambda$, and let
$\wbar_{\lambda}( X)$ be its closure.  For example, $w_{1^n}(X)$ is
the configuration space of $n$ unordered distinct geometric points, sometimes
denoted $B(X,n)$ or $\Conf^n(X)$.  Let $1^k \lambda$ denote the partition obtained
from adding $k$ $1$'s to $\lambda$.

\tpoint{Theorem}\label{T:IntroPointsLimits}
{\em 
Suppose $X$ is a geometrically irreducible variety of dimension $d$.
\begin{enumerate}
\item[(a)] If $X$  satisfies MSSP${}_\phi$, then the limits
\begin{equation}\label{eq:AA}
\lim_{j\ra\infty}  \frac{[w_{1^j\lambda} (X)]}{\L^{dj}}
\qquad
\textrm{and}
\qquad
\lim_{j\ra\infty}  \frac{[\wbar_{1^j\lambda} (X)]}{\L^{dj}}
\end{equation}
   exist in $\phi(\K)$,
and have finite formulas in terms of motivic zeta values, the $[ \Sym^i
X]$, and 
$\Sinf_\phi(X) $ (defined in \S \ref{c:mssp}).
 If furthermore  the $[\Sym^j X]$ are invertible in $\phi(\K)$
(e.g.\ if $X$ is rational or $\phi = HS$, \S \ref{o:Syminv}), then 
\begin{equation}\label{eq:BB}
\lim_{j\ra\infty}  \frac{[w_{1^j\lambda} (X)]}{\left[ \Sym^{j+|\la|} X \right]}
\qquad
\textrm{and}
\qquad
\lim_{j\ra\infty}  \frac{[\wbar_{1^j\lambda} (X)]}{\left[
\Sym^{j+|\la|} X \right]}
\end{equation}
   exist in $\phi(\K)$,
and have finite formulas in terms of motivic zeta values and $[ \Sym^i
X]$.
\item[(b)]  If $\k = \F_q$, then 
\begin{equation}\label{eq:CC}
\lim_{j\ra\infty}  \frac{\# w_{1^j\lambda} (X)}{ \#\Sym^{j+|\la|} X }
\qquad
\textrm{and}
\qquad
\lim_{j\ra\infty}  \frac{ \# \wbar_{1^j\lambda} (X)}{ \#
\Sym^{j+|\la|} X }
\end{equation}
   exist,
and have finite formulas in terms of Weil zeta values, and derivatives
of the Weil zeta function.
\end{enumerate}
}

The second statement of (a) can be interpreted as a limiting ``motivic
probability'' ($w_{1^j \la}(X)$ and $\wbar_{1^j \la}(X)$ are both
subsets of $\Sym^{j + |\la|} X$), and the statement of (b) can be
interpreted as a limiting probability.  Corollary~\ref{C:stabform}
is a more precise version of Theorem~\ref{T:IntroPointsLimits}.
Corollary~\ref{C:stabform} in turn is a consequence of an analogous
statement about generating functions (Theorem~\ref{T:Kform}), with
{\em no}
motivic stabilization hypotheses, which states that $\sum_{j} [w_{1^j \lambda} (X)]
t^j$ and $\sum_{j} [\wbar_{1^j \lambda} (X)] t^j$ have finite formulas
in terms of $\zzeta_X(t)$ and $[ \Sym^i X]$; the  formulas are essentially
the same as in Theorem~\ref{T:IntroPointsLimits} above.

The finite formulas of Theorem~\ref{T:IntroPointsLimits} are given
recursively in Propositions~\ref{P:Krecursion} and~\ref{P:K}(b).  It
is not hard to show that the limits exist; the main content of
Theorem~\ref{T:IntroPointsLimits} (or, rather, Corollary~\ref{C:stabform}) is the description of the limit.
These limits have particularly nice descriptions in special cases.  We
give some now.  Rather than giving three forms of each result
(corresponding to \eqref{eq:AA}--\eqref{eq:CC} in 
Theorem~\ref{T:IntroPointsLimits}), for simplicity we just discuss the
``motivic probability'' versions (\eqref{eq:BB} in Theorem~\ref{T:IntroPointsLimits}), as representative of all three
versions.

\point \label{pre56}Fix now a rational variety $X$ (over $\k$) of
dimension $d$.   (The reason for assuming rationality is for simplicity, so that $[\Sym^j X]$ is invertible in $\K$, \S
\ref{o:Syminv}.)   The limiting motivic probability (as $j \rightarrow \infty$)
that $j$ points are distinct (i.e.\ $\lim_{j \rightarrow \infty} [
w_{1^j}(X) ] / [ \Sym^j X]$) is $ \bzeta_X(2d)^{-1}$.  (Equivalently, the
limiting motivic probability that $j$ points are not distinct --- the
traditional ``discriminant locus'' --- is $1 -
\bzeta_X(2d)^{-1}$.)    The corresponding generating function formula 
is
$$
\sum_{j} [w_{1^j} (X)] t^j =\zzeta_X(t)/\zzeta_X(t^2)
$$
(a special case of Proposition~\ref{P:Kbase}(a)), which specializes,
under taking Euler characteristic with compact supports, to the well-known formula for Euler characteristic of configuration spaces
$$
\sum_{j} \chi_c(w_{1^j} (X)) t^j =(1+t)^{\chi(X)}
$$
using Macdonald's formula $\chi_c(\zzeta_X(t))=(1-t)^{-\chi(X)}$ from
\cite{Mac62}.  We remark that there is a large body of work, going
back to Macdonald \cite{Mac62}, giving generating functions for
motivic or topological invariants of symmetric products (see
\cite{Z72, M78, Cheah96, BL02,O08, MS11, MS12}) and Hilbert schemes
(see \cite{G90, G93, Cheah96, BL03,GZLMH04, GZLMH06, BNW07, NW08,
  CMOSY12}).  Our formulas also extend such motivic generating
functions to the generalized configuration spaces $w_{1^j\lambda} (X)$
and $\wbar_{1^j\lambda} (X)$ which are the natural strata (and their
closures) of symmetric products.

\point \label{promise56}We return to our examples.  Generalizing \S
\ref{pre56}, the limiting motivic probability that $j$
unordered points have a point of multiplicity (at least) $a$, i.e.\
$\lim_{j \rightarrow \infty} [ \wbar_{1^{j-a}a} (X) ] / [\Sym^j X]$, is $1
-  \bzeta_X(ad)^{-1}$ (a consequence of Proposition~\ref{P:Kbase} and
Lemma~\ref{L:limits}). 
  This is an analog of the classical arithmetic fact that the
proportion of $a$th-power-free integers is $\zeta(a)^{-1}$.
% (and is
%similarly straightforward).

\point \label{promise513} More generally (and more subtly) there is a simple
description of the limiting motivic probability that there are $r$ points of
multiplicity $b$ or worse, i.e.\ $\lim_{j \rightarrow \infty} [
\wbar_{1^j b^r} ] / [\Sym^{j+br} X]$ (a consequence of
Proposition~\ref{P:Kjksrecursion} and Lemma~\ref{L:limits}); its simplicity is
clearest in the case where $X = \A^d$, in which case $[\wbar_{1^j b^r}(\A^d)]/[\Sym^{j+br}\A^d]=1 / \L^{d
  r(b-1)}$  for $j\geq 0$ (see Example~\ref{ex:last}).
As an even more specific example, the probability that a polynomial
(of degree at least $4$) over $\FF_q$ has two double roots ``or
worse'' (a quadruple root; a triple root is {\em not} enough)
is $q^{-2}$.

\point \label{promise59}As a further example, if $\nu$ has all distinct elements greater than
$1$, then
$$\lim_{j\ra \infty}\frac{ [ w_{1^j\nu} (X)] }{[ \Sym^{j+\sum \nu} X ]} = 
\frac { [ w_\nu  (X) ] } {\bzeta_X (2d)}
 \frac{   \L^{-d \sum
    \nu}}{ \left(1+\L^{-d}\right)^{|\nu|}}.$$
 (See
Example~\ref{E:distinctnu}.)
Using the ``fibration''  $\alpha: w_{1^j \nu} \rightarrow w_\nu$, one can give a 
``fiberwise heuristic'' which yields this as a prediction.
 But because $\alpha$ is not a fibration in the Zariski topology, this
heuristic does not give a proof, so we omit the details.

\point Our last specific example is the following.
\label{d:syms}Let \highlight{$\Sym^{j}_{s} X$} (not to be confused with
$\Sym^j_{[s]} X$, \S \ref{s:mzf}) be the locally closed subset of $\Sym^j X$
corresponding to collections of points containing {\em exactly} $s$
multiple points.

\tpoint{Theorem}
{\em \label{T:IntroPointsssing}
Suppose $X$ is a geometrically irreducible variety of dimension $d$.
\begin{enumerate}
\item[(a)]
If $X$ satisfies MSSP${}_\phi$,
  then
$$
\lim_{j\ra\infty} \frac{ \left[ \Sym^{j}_{s} X \right]}{\L^{j d}}=
\frac { \bzeta^{[s]}_X({2 d} )}   {\bzeta_X({2 d})}\Sinf_\phi
(X) \qquad\text{    in $\phi(\K)$}.
$$
 If furthermore  the $[\Sym^j X]$ are invertible in $\phi(\K)$
(e.g.\ if $X$ is rational or $\phi = HS$, \S \ref{o:Syminv}), then 
\begin{equation}\label{eq:blimp}
\lim_{j\ra\infty} \frac{\left[ \Sym^{j}_{s} X \right]}{ \left[ \Sym^j
    X \right]}=
\frac { \bzeta^{[s]}_X({2 d} )}   {\bzeta_X({2 d})}  \qquad\text{    in $\phi(\K)$.}
\end{equation}
\item[(b)] If $\k = \F_q$, 
$$\lim_{j\ra\infty} \frac{ \# \Sym^{j}_{s} X }{  \# \Sym^j
    X }=
\frac { \zeta^{[s]}_X({2 d} )}   {\zeta_X({2 d})} .
$$
\end{enumerate}
}
The proof of Theorem~\ref{T:IntroPointsssing} concludes just after the
statement of Theorem~\ref{T:Kssing}.  

The similarity of \eqref{eq:blimp} to Theorem~\ref{T:Bjorn} is striking.
In fact, for $X=\PP^1$, Theorem~\ref{T:Bjorn} and
Theorem~\ref{T:IntroPointsssing} give the same result, but for a
smooth  curve $C$ of arbitrary genus, Theorem~\ref{T:IntroPointsssing} gives the limit
of the moduli spaces of \emph{all} divisors, and Theorem~\ref{T:Bjorn}
gives the analogous result for moduli spaces of divisors in multiples
of a fixed linear system (although the answer does not depend on the
linear system). 
Thus although in the case of (smooth projective  geometrically
irreducible) curves of positive genus,
 \eqref{eq:blimp} and Theorem~\ref{T:Bjorn} are logically independent,
 they are consistent in some strong sense.

\bpoint{Connections to configuration spaces in topology}\label{S:top}

We now draw connections to topological work.  

\epoint{The ``contractible'' case $X =\A^d$}\label{S:contract} We begin with the case
where $X= \A^d$, to highlight the topology
arising from the positions of the points rather than the underlying
space.  Note first that $[\Sym^r \A^d] = [\A^{rd}]$ (a fact first
proved by Totaro, \cite[Lemma 4.4]{Got01}, see also \cite[Thm.~1]{GZLMH11}
 and \cite[Statement~2]{GZLMH04}), so $\zzeta_X(t) = 1 / (1- [X]
t)  = 1 / (1- \L^d
t)$.

Proposition~\ref{P:Kbase}(b) implies that $[ \wbar_{1^j a} (\A^d)] =
\L^{d(j+1)}$.
Occam's
Razor~\ref{s:occraz} then gives a striking prediction (the case $r=0$ of
Conjecture~\ref{c:q} below).
The results given in \S
\ref{promise513}, and more generally Example~\ref{ex:last}, suggest
even more (the full statement of Conjecture~\ref{c:q} below).

  Before stating it, we point out that for an arbitrary complex manifold
$X$, it is more natural to study the complement $\wbar^c_\lambda(X)$
of $\wbar_\lambda(X) $ in $\Sym^{\sum_i \lambda_i} (X)$.  For example,
$\wbar^c_{1^j 2} (X)=w_{1^{j+2}} (X)$.  The spaces $\wbar^c_\lambda(
X)$ satisfy Poincare duality (with $\QQ$-coefficients, see Motivation~\ref{motivation}(v)).

\tpoint{Conjecture~E}
{\em If $1 < a \leq b$, and $j$ and $r$ are nonnegative integers,
% then
% $h^i(\wbar_{1^j a b^r}(\A_{\C}^d), \Q)=0$ for $i>0$. 
then $h_i(\wbar^c_{1^j a b^r}(\A_{\C}^d), \Q)$ is $1$
  if $i=0$ or $i=2d((a-1) + r (b-1))-1$, and $0$
  otherwise.}\label{c:q}

\point \label{s:arnold}In the case when $d=1$, $a=2$, $r=0$ (hence $b$ arbitrary), and $j$
arbitrary, this is a result of Arnol'd
\cite{Arnold69}.  In the case when $d=1$, $a=b$, and $r$ and $j$ are
arbitrary, this is a consequence of \cite[(19) and (20)]{a2} (but note
a mistake in the formulation of \cite[(23)]{a2}).  
T. Church explained this to us, and explained how Arnol'd's
proofs of these cases extends to general $d$, \cite{churchpc}. 
O. Randal-Williams
\cite{RWpc}
has also proved the case of one multiple point (i.e.\  $r=0$,
and $j$, $a$, $b$,  and $d$ arbitrary) using \cite{sk}.% Soren also did

Given Example~\ref{ex:last}, the reader may suspect that $\wbar_{\la}(\A^d)$
is always a power of $\L$ for all $\la$, but this is
not the case.  The smallest $\la$ for which this is false is $\la = 1^2
2^2 3$ (see the last line of \S \ref{s:genpar}).  In fact, as $j
\rightarrow \infty$, $\wbar_{1^j 2^2 3}(\A^d_{\C})$ has an unbounded
number of nonzero cohomology groups with compact support; this can be seen through a
calculation of the generating series $\sum_{j} [\wbar_{1^j 2^2 3}
(X)] t^j$ using Theorem~\ref{T:Kform}.

\epoint{General $X$}
The stabilization of the Betti numbers (in fact the integral homology) of $w_{1^j}(X)$
 for open manifolds $X$ was proven by McDuff \cite{McDuff75}.  
Recently,  Church \cite[Cor.~3]{Church12} and Randal-Williams \cite{RW12} proved the 
stabilization of the Betti numbers $h_k(w_{1^j}( X),\QQ)$ for closed,
connected manifolds $X$ of finite type.
This is the topological analog of the motivic limits existing in our
first example: $\lim_{j\ra\infty}  [w_{1^j} (X)] / \L^{jd}$ (cf.\
Occam's Razor~\ref{s:occraz}).
Upon hearing of our result for the motivic stabilization of partially
labeled configuration spaces $[w_{1^j\lambda} (X)] / \L^{-jd}$, Church \cite[Thm.~5]{Church12}  and Randal-Williams \cite{RWpc} also proved the stabilization of the Betti numbers of these spaces for manifolds $X$.

We conjecture stabilization of Betti numbers for the other flavors of
configuration spaces that have motivic limits (Theorems~\ref{T:IntroPointsLimits} and \ref{T:IntroPointsssing}).

\tpoint{Conjecture~F} {\em 
Given $i$ and a partition $\lambda$, for an irreducible smooth complex variety $X$, the limit
  $\lim_{j\ra\infty} h_i( \wbar^c_{1^j \lambda}(X), \Q)$ exists. } % Oscar can do $r=1$

See \S \ref{s:arnold} for the case $X =
\A^d_{\C}$ and $\lambda = m^r$ ($m$ and $r$ arbitrary).
 (As in the case of $w_{1^j} (X)$, there do not exist
obvious maps among the elements of these sequences of configuration
spaces for closed $X$; many topological stabilization results rely on
such a map.)  
Since the authors publicly made this conjecture in the first draft of this paper,
Kupers and Miller have proven Conjecture F in the case that
$\lambda=m$ for an integer $m$, for $X$ a ``reasonable'' manifold
\cite{KM2}.  This can be done using the topological methods of
\cite{McDuff75}, but Kupers and Miller  can also show homological stability
with a specific range, and moreover their result holds with
$\Z$-coefficients so long as the manifold is not closed.  (The
case when $\dim X=2$ and $\lambda=m$ was done earlier in \cite{y}.  We also note
that Baryshnikov has studied  the topology of the space
$\wbar^c_{1^jm}(X)$ with a view towards applications, see \cite{bar}.)  Moreover, the
methods
of Kupers and Miller may apply more generally to partitions of the
form $\la = m^r$.

One might ask a similar question for the constructible
subset $\Sym^j_s X \subset \Sym^j X$, for each $s$ (cf.\ Theorem~\ref{T:IntroPointsssing}).

The formulas of Theorem~\ref{T:IntroPointsLimits} (given recursively
in Propositions~\ref{P:Krecursion} and~\ref{P:K}(b)) can be combined
with the formulas for the Hodge-Deligne series of zeta functions
\cite[Prop.~1.1]{Cheah96} to obtain explicit formulas for the
Hodge-Deligne series of the limits of (various flavors of)
configuration spaces above ($w_{1^j \lambda}(X), \wbar^c_{1^j \lambda}(X)$, and $\Sym^{j}_{s} X$ ) in terms of the Hodge-Deligne polynomial of
$X$.  Totaro \cite{T96} gives an explicit spectral sequence with only
one non-trivial differential that computes the Betti numbers of the
usual configuration spaces, but this does not immediately give the
limit Poincar\'e series.  One should hope to compare the limit
Hodge-Deligne series of various configuration spaces to the analogous
(mostly unknown) limit Poincar\'e series.  In particular, in the
situations in this paper, we not only know that limit Hodge-Deligne
series exist, but we have given relatively simple formulas for them.
Are there analogous simple formulas for the limit Poincar\'e series?

More precisely, in analogy with Conjectures~C and~D, for all the
flavors of configuration spaces we discuss, we conjecture the limiting
Poincar\'e series is finite when the analogous limiting motive has
finite Hodge-Deligne series, and wonder whether, in these cases,
Occam's Razor~\ref{s:occraz} predicts the correct Betti numbers. 
 As
stated in \S \ref{s:arnold}, in the case $X=\A^d_\C$, for
$\wbar^c_{1^j 2}$,
Arnol'd has shown the answer is
yes; and for $X=\P^1_\C$, for $w_{1^j}$, Church has done the same, \cite[Prop.~4.5]{Church12v1}.

As an  example, applying Occam's Razor~\ref{s:occraz} to
Example~\ref{pre56} in the case $X=\P^2_\C$ yields the following prediction.

 \tpoint{Conjecture~G}{\em 
We have 
$$
\lim_{j\ra\infty} h_i(w_{1^j} (\P^2_\C),\QQ)=
\begin{cases}
1 &\textrm{ if $i=0,2,4,7,9,11$} \\
0 &\textrm{ otherwise.} 
\end{cases}
$$
}

Since the authors publicly made this conjecture in the first draft of this paper,
Kupers and Miller \cite{KM} have proven Conjecture G using the scanning map of McDuff \cite{McDuff75}. 
However, Conjecture G is just one of many possible cases that illustrate our conjecture that
the limiting
Poincar\'e series is finite when the analogous limiting motive has
finite Hodge-Deligne series.  For example, we also have the following from applying Occam's Razor~\ref{s:occraz} to
Example~\ref{promise59} in the case $X=\P^1_\C$ and $\nu=23$.
 \tpoint{Conjecture~G'}{\em 
We have 
$$
\lim_{j\ra\infty} h_i(w_{1^j23} (\P^1_\C),\QQ)=
\begin{cases}
1 &\textrm{ if $i=0,1$} \\
0 &\textrm{ otherwise.} 
\end{cases}
$$
}
While any particular case might be handled by existing topological methods, we are curious whether there are topological methods that could handle all or many cases at once.

\point We conjecture
that the limiting Poincar\'e series is periodic when
the analogous limiting motive has periodic Hodge-Deligne series, as in
the following example, which is about the space of configurations with precisely
one double point.
 \tpoint{Conjecture~H} {\em
The limits $\lim_{j \ra\infty} h_{i}(w_{1^j2}(\A^d_\C),\QQ )$ are periodic in $i$.
}

If we further apply Occam's Razor~\ref{s:occraz} to our results, it would predict that for each $i$,
\begin{equation}\label{badoccam}
\lim_{j \ra\infty}  h_{i}(w_{1^j2}(\A^d_\C),\QQ ) =
\begin{cases}
 1 &\textrm{if $i=0$}\\
 2 &\textrm{if $i=2(2k-1)d -1$ or  $4kd$, for $k\geq 1$}\\
 0 &\textrm{otherwise.}
\end{cases}
\end{equation}
%O. Randal-Williams has suggested that Conjecture~H is true but that
%the prediction \eqref{badoccam} is false, \cite{RWpc}.
Since the authors publicly made this conjecture in the first draft of this paper,
Kupers and Miller \cite{KM} have proven Conjecture H and shown that the prediction \eqref{badoccam} is false by finding the actual stable Betti numbers.

Tommasi has made the following beautiful observation since the first draft of this paper appeared.  It is known that the rational cohomology
of $w_\lambda(\A^d_\C)$ is non-trivial only in degrees that are a multiple of $2d-1$ and 
 that the cohomology groups in degree $k(2d-1)$ are Hodge structures of Tate of weight $2kd$. 
(This can be seen when the elements of $\lambda$ are distinct from \cite{T96}, and then this fact passes to finite quotients.)  So in fact for $w_\lambda(\A^d_\C)$ the Hodge-Deligne polynomial determines the Poincar\'e polynomial.  More precisely, we have that 
$$
\left.\frac{e(w_\lambda(\A^d_\C))}{(xy)^{d|\lambda|}}\right|_{(xy)^{-d}=-t^{2d-1}}=\sum_k
h_k(w_\lambda(\A^d_\C)),\Q) t^k.
$$  
(The quotient by $(xy)^{d|\lambda|}$ is from using Poincar\'e duality to pass from cohomology with compact supports to cohomology.)

So, Tommasi's observation plus our Example~\ref{E:distinctnu} imply Conjecture H and in fact
\begin{equation}\label{goodoccam}
 h_{i}(w_{1^j2}(\A^d_\C),\QQ ) =
\begin{cases}
 1 &\textrm{if $i=0,j(2d-1)$}\\
 2 &\textrm{if $i=k(2d-1)$ , for $1<k<j$}\\
 0 &\textrm{otherwise.}
\end{cases}
\end{equation}
Moreover, Tommasi's observation implies  our conjecture that the limiting
Poincar\'e series is finite when the analogous limiting motive has
finite Hodge-Deligne series (given before Conjecture G)  in the cases $w_\lambda(\A^d_\C)$, and
while the stable Betti numbers are not given by Occam's Razor~\ref{s:occraz}, they are given by a simple transformation of the Hodge-Deligne polynomials for which we determine a generating function.  Further, Tommasi's observation implies our conjecture
that the limiting Poincar\'e series is periodic when
the analogous limiting motive has periodic Hodge-Deligne series (given before Conjecture H) in the cases $w_\lambda(\A^d_\C)$.  
However, there are still many similar patterns that are predicted by our results on the limiting Hodge-Deligne series, such as the following.

 \tpoint{Conjecture~H'} {\em
The $i$ for which the limits $\lim_{j \ra\infty} h_{i}(w_{1^j22}(\P^1_\C),\QQ )$ are non-zero are periodic, and the non-zero limits are $1$.
}

\bpoint{Connections to configuration spaces in number theory}\label{S:arith}

The limits in Theorems~\ref{T:IntroPointsLimits} and \ref{T:IntroPointsssing}
have natural analogs over $\Z$.  
For a partition $\nu=[e_1,e_2,\dots,e_k]$, we say an integer $n$ has \emph{at least $\nu$-power} if 
$\prod_i a_i^{e_i} | n$ for some (not necessarily distinct) integers $a_i>1$.
The limit resulting from the generating function in 
Example~\ref{ex:last} then has the following analog over $\Z$:
\begin{align*}
 &\lim_{N\ra\infty} \frac{\#\{1\leq n\leq N\; |\;  \textrm{$n$ has at least $ab^r$-power}\}}{N}\\&=
1- \frac 1 {\zeta(-b)}\left(\sum_{i=0}^{r-1} 
\sum_{p_1\leq\dots\leq p_i} p_1^{-b}\cdots p_i^{-b}\right)
-\sum_{p_1\leq\dots\leq p_r} p_1^{-b}\cdots p_r^{-b} \frac 1 {\zeta(-a)},
\end{align*}
where the sums above are over primes $p_j$.
Theorems~\ref{T:IntroPointsLimits} and \ref{T:IntroPointsssing} also
 suggest natural point counting analogs  for arithmetic schemes 
(as in  \cite[\S 5]{Poonen04}).  
One expects that when $X$ is a general arithmetic scheme, 
such results, as in \cite[\S 5]{Poonen04}, will require new ideas.

\bpoint{Acknowledgements} 

We have benefited from conversations from
many mathematicians.  
In particular, we would like to thank the following for helpful conversations 
and comments on earlier versions of this paper: 
D. Arapura,
T.  Church, R. Cohen, J. Ellenberg, D. Erman, S. Galatius, E. Howe, A. Kupers,
C. Liedtke, J. Miller, B. Poonen,  
O. Randal-Williams,   O. Tommasi, B. Totaro,  A. Venkatesh,  and C. Westerland.
\begin{comment}
  T.  Church, R. Cohen, S. Galatius, O. Randal-Williams,
and C. Westerland patiently explained facts from topology to us.
J. Ellenberg, E. Howe, B. Poonen, and A. Venkatesh did the same for
arithmetic geometry.  D. Arapura, D. Erman, and C. Liedtke did the
same for algebraic geometry.  B. Totaro did the same in several
fields.  We thank B. Poonen for the proof of Lemma~\ref{L:SymInd}.
Ellenberg, Venkatesh, and Westerland's paper \cite{evw} is a wonderful
example of a connection between the style of questions in topology and
style of questions in arithmetic that our motivic work has connections
to on each side.  B. Poonen and O. Randal-Williams gave helpful
comments on an earlier version of the paper.  
\end{comment}
The authors gratefully
acknowledge the support of NSF grant DMS-1100771 (RV), and an AIM
Five-Year Fellowship and NSF grant DMS-1001083 (MMW).

\section{Notation for partitions and configuration
  spaces} 
\label{s:genpar} 

Our arguments use partitions in a slightly more general sense than usual.  For
us, a \highlight{\emph{(generalized) partition}} in an abelian
semigroup $S$ is a finite multiset of elements of $S$.  A {\highlight{\em
  subpartition}} is a submultiset of a partition.  Partitions in the
traditional sense are the special case $S = \Z^+_{>0}$.  
We use the standard notation \highlight{$\mu \vdash n$} (``$\mu$ is a
partition of $n$'').
Let
\highlight{$\cP$} be the set of all partitions in the traditional
sense (i.e. of positive integers). 
We use the notation \highlight{$[ \cdots ]$} to denote a multiset,
exponents to denote multiplicity, and concatenation to denote union,
so for example $a^2 b = a [a,b] = [a] [a, b] = [ a,a,b ]$.

In \S \ref{S:Hyper}, we will want to concatenate 
$\la$ and $\mu$ and consider the parts of $\la$ and $\mu$ to
be distinct (``disjoint concatenation''), which may require some renaming; we write this as $\la
\conc \mu$.  For example, we may write
$[a,b] \conc [a,a]$ as $[a_1, b_1, a_2, a_2]$. 

As usual,
\highlight{$|\l|$} denotes the number of elements of the multiset $\l$.
A generalized partition $\lambda$ has a \highlight{{\em multiplicity partition} $m(\lambda)$} of $|\l|$ --- for example $m( [  a, a, b ] )= [  2, 1 ]$.  
We write \highlight{$||\lambda||$} for the number of {\em distinct} elements of
$\lambda$. %(the number of $i \in S$ for which $m_i\ne 0$).
We write \highlight{$\sum \lambda$} for \emph{the sum} of the
generalized partition $\lambda$, i.e.\ $\sum_{s\in \lambda} s$.  
Clearly, for any generalized partition, $|\l| = \sum m(\l)$ and $||
\l || = | m(\l) |$. For example, 
suppose $\lambda=1^3 2^3 3 4^2 5$, so $m(\l) =[3,3,1,2,1 ]$.  Then $\sum
\lambda=25$,  $|\lambda|=\sum m(\l)=10$,  $||\lambda||=|m(\l)|=5$, and $||m(\l)||=3$.

Let \highlight{$\cQ(m)$} be the set of partitions in the
traditional sense in which exactly the numbers $1$ through $m$ appear 
(so for example
$||\mu|| = m$ for $\mu \in \cQ(m)$).
(We also also think of $\cQ(m)$ as partitions of $m$ linearly ordered elements, up to isomorphisms of the ordered elements.)
By taking the multiplicity partitions, we can interpret $\cQ(m)$ as the {\em ordered} partitions with exactly $m$
parts.  For example, $[1,1,1,1,2,3,3] = 1^4 2^1 3^2\in \cQ(3)$ can be reinterpreted
as $4+1+2=7$.  Let \highlight{$\cQ = \cup_m \cQ(m)$}, which can be reinterpreted
as the set of all ordered partitions.

Suppose $\lambda$ and $\lambda'$ are generalized partitions in $S$.
If there are sub(multi)sets $[x,y] \sub\lambda$ and $[z]\sub\lambda'$ such
that $x+y=z$ and $\lambda\setminus [x,y]= \lambda'\setminus[z]$, we
say $\lambda'$ is an \highlight{{\em elementary merge}} of $\lambda$.
In this case $|\lambda|=1+|\lambda'|$. 
We define the \highlight{{\em refinement ordering}} $<$ on generalized
partitions  in $S$ as generated by elementary merges.  (If $\lambda'$ is
an elementary merge of $\lambda$, then 
\highlight{$\lambda < \lambda'$}.)
 For example, $[  1,2,3 ] < [  3, 3 ] < [  6 ]$.  We write  \highlight{$\la \leq \la'$} if $\la < \la'$
or $\la  = \la'$.

Given a generalized partition $\lambda=[\la_i]$, define the \highlight{\emph{formalization}}
as \highlight{$f(\lambda) :=[a_{\la_i}]$} (in the abelian semigroup
$\Z^+[a_i]_{i \in S}$); we have replaced entries with ``formal''
replacements.   The purpose of this construction is to obtain a partition
with the same multiplicity sequence such that 
for all $S_1,S_2\subset \lambda$ such that $\sum S_1=\sum S_2$, we
have $S_1=S_2$.
   %(This can be defined in terms of a universal property.)
%Note that $[  1, 1, 2, 3, 3, 3 ]$ is not formal (as $1+1=2$), but
%its formalization $[  a_1, a_1, a_2, a_3, a_3, a_3 ]$ {\em is} formal (note that $a_1 + a_1 \neq a_2$).
%We define the  \highlight{\emph{formal  merge (partial) ordering $\ll$}} as
%generated by 
%\highlight{$\lambda \ll  \lambda'$}
%if
%$f(\lambda)<\lambda'$.
%For example,
%$$
%[1,1,2,3] \ll  [2a_1, a_2, a_3] \ll [a_{2 a_1 + a_3}, a_{a_2}].
%$$

If $\l$ is a generalized partition, define \highlight{$\Sym^{\la} X = \prod_{m_i \in m(\l)} \Sym^{m_i}
   X$}.  
 For
 example, $\Sym^{ [  a, a, b ]  }X$ parametrizes an unordered pair of
 (geometric, not necessarily distinct) points of $X$ labeled $a$, and another point (not
 necessarily distinct)   labeled $b$.
(Warning: do not confuse $\Sym^{[2]} X$ with $\Sym^2 X$: by definition
$\Sym^{[2]}X = X$.)    We 
define  $w_\l( X)$ (or simply $w_\l$ for convenience) 
to be  the open  subscheme of
$\Sym^{\l}X$ in which all the points are distinct, i.e.\ the
complement of the ``big diagonal''. 
(This generalizes the definition of $w_{\l}$ given at the start of \S
\ref{s:conspa}, which is  the case of traditional  partitions.)
 For example, $w_{[  aab ] }(X)$
parametrizes an unordered pair of of distinct points of $X$  labeled
$a$, along with a third distinct point, labeled $b$.   
Note that $w_\l$  depends only on the multiplicity sequence $m(\lambda)$.

Define  \highlight{$\overline{w}_\lambda = \sum_{\lambda \leq \mu}
  [w_\mu]$}.
Although $\wbar_\l$ is defined as an element of $\cM$, we
can often naturally endow it with the structure of a variety, as the
closure of $w_\l$ in an appropriate space.  For example, if $\l$ is a
traditional partition ($S= \Z^+$), then $\wbar_\l$ is the class of the closure of
$w_\l$ in $\Sym^{\sum \l} X$; thus this definition of $\wbar_\l$
generalizes the one given  at the
start of \S \ref{s:conspa}. 
The varieties $w_\la$ and $\wbar_{\la}$ have been studied by Haiman and Woo
(see $Z^\circ_\lam$ and $Z_\la$ in \cite[\S 3.2]{HWpre2012}).
 If $\l$ is a formalization,  since
$w_{\lambda}$ is  the open subset of
$\Sym^{\lambda} X$ where the $| \lambda |$ points are distinct, and
the various $\mu$ with $\lambda < \mu$ correspond to letting the points
come together in various ways, we have
\begin{equation}
\overline{w}_\lambda(X) =\left[ \Sym^{\lambda}X \right] \quad  \textrm{ for  a formalization $\l$.} \label{e:formal}
\end{equation}
But \eqref{e:formal} need not hold if $\l$ is not a formalization.  As
perhaps the simplest example, if $\l = [1,1,2,2,3]$, then
$\wbar_{\l}(\L) = \L^5 - \L^2 + \L$, while $\Sym^{\l} \L = \L^5$ (using $\Sym^n \L = \L^n$).
\hidden{ (Side remark: $\wbar$ can be endowed with an interpretation
  as a variety even if $\l$ is not formal or does not consist only of
  integers.)  }

\section{Moduli of hypersurfaces}\label{S:Hyper}

 The goal of this section is to prove Theorem~\ref{T:Bjorn}.
Throughout this section, $X$ is assumed to be smooth of pure dimension
$d$. 
In order to prove Theorem~\ref{T:Bjorn} in general, we first establish
it for $s=0$.  This case  will be completed by 
Proposition~\ref{P:homog}, see \S \ref{s:s0done}.   We determine the motive of smooth divisors in a
linear system by considering all divisors, and removing those with
singularities.

Suppose \highlight{$\l$} is a generalized partition and \highlight{$\cF$} is a line bundle on
$X$.
We define three types of incidence subschemes parametrizing
sections of $\cF$ singular  at points marked by $\l$.

In analogy with the notation $w_\l$, let \highlight{$W_\l(\cF)$} (or $W_\l$ when $\cF$ is clear from context) denote the
locally closed subvariety of $ H^0(X,\cF) \times w_\l(X) $ corresponding
to sections of $\cF$ singular at {\em precisely} those $|\l|$
(necessarily distinct) geometric points of $X$ given by the point of $w_\l(X)
\subset \Sym^{\l}X$.   For example, $W_{*^s}(\cF) \cong H^0(X, \cF)^s$.

Let \highlight{$W_{\geq \l} = W_{\geq \l}(\cF)$} be the locally closed
subvariety of $ H^0(X,\cF) \times w_\l(X)$ corresponding to sections of $\cF$
singular at those $|\l|$ (necessarily distinct) geometric points of $X$ given by
the point of $w_\l(X) \subset \Sym^{\l}X$, {\em and possibly elsewhere}.
Note that $W_{\la}$ is an open subset of $W_{\geq \la}$.

If $k$ is a nonnegative integer, let \highlight{$W_{\lambda,{\geq k}}= W_{\lambda,{\geq k}}(\cF)$} be
the locally closed subset of $ H^0(X, \cF) \times w_\l(X)$
corresponding to those $(s,t)$ for which $s$ is singular at the
$|\la|$ points
parametrized by $t$ and {\em at least $k$ additional geometric  points}.
Because $W_{\la, \geq k}$ is the image of $W_{\geq \la \cdot *^k}$
(disjoint concatenation
``$\conc$''  was defined in \S \ref{s:genpar}) under the
obvious projection, $W_{\la, \geq k}$ is a constructible subset of 
$ H^0(X, \cF) \times w_\l(X)$ by Chevalley's Theorem, and thus
has a well-defined class in $\cM$.
(This also follows
 from \eqref{E:breakintoZs} below.)

Clearly
\begin{equation}\label{E:breakintoZs}
 [W_{\geq \lambda}] = [W_{\l}]+ [W_{\lambda,{\geq 1}} ]
= [W_{\l} ] +  [  W_{\l \conc *}  ] + [  W_{\lambda,{\geq 2}}]
= [W_{\l} ] + [   W_{\l \conc *} ] + [ W_{\l \conc **}  ] + [
W_{\lambda,{\geq 3}}  ] = \cdots.
\end{equation}
For example, a section singular at some points labeled by $\l$ is:  
(0) nonsingular elsewhere, or else (i) singular at precisely one point
elsewhere, or else (ii) singular
at precisely two points elsewhere, or else (iii) singular at $3$ {\em or more}
other points elsewhere.

\tpoint{Lemma} {\em With $\cL$ ample and fixed, and $j$ sufficiently
  large in terms of $|\l|$, we have that $W_{\geq \lambda}(\cL^{\otimes
    j})$ is a vector bundle over $w_\l(X)$ of rank $r-|\l|(d+1)$,
  where \highlight{$r = h^0(X,\cL^{\otimes j})$}.\label{l:bundle} }
 
\noindent {\em Proof.}  The following argument will not surprise experts, but we include it for completeness.
The result is insensitive to base field extension, so we assume  $\k = \overline{\k}$.  
The scheme $W_{\geq \l}$ corresponds to a coherent sheaf on $w_\l(X)$, corresponding to sections of   $\cL^{\otimes j}$
singular at the $|\l|$ points parametrized by $\Sym^{\l} X$.  We wish
to show that this coherent sheaf is a vector bundle of rank $h^0(X,
\cL^{\otimes j}) - |\l| (d+1)$.

By Grauert's Theorem, it suffices to show that (for $j \gg_{|\l|} 0$)
for any closed point of $w_\l(X)$, interpreted as $|\l|$ distinct
points of $X$, the $1$-jets at the points impose independent
conditions on sections of $\cL^{\otimes j}$.  Now choose $j$ so
$\cL^{\otimes j}$ is $(|\l|(d+1))$-very ample.  
\qed
 
The typical Noetherian induction using the local triviality of vector
bundles yields, 
 for $j \gg_{|\l|} 0$,
\begin{equation}\label{eq:W}
[W_{\geq \lambda}] =[w_\l ] {\mathbb L}^{r-|\l|(d+1)}. % = w_\l ( \L^{r-|\l|(n+1)}  + \cdots + 1).
\end{equation}

\tpoint{Corollary (and definition of $j_N$)}
{\em Fix an ample line bundle  $\cL$ on $X$ of dimension $d>0$.
For each positive integer $N$, there is some \highlight{$j_N$}, so that for $j
\geq j_N$ (where as in Lemma~\ref{l:bundle}, $r = h^0(X,\cL^{\otimes j})$):\label{c}
\begin{enumerate}
\item[(a)]  $r>2N$; 
\item[(b)] $W_{\geq \lambda}(\cL^{\otimes j})$  
(and hence its open subset $W_\l(\cL^{\otimes j})$)
has (pure) dimension $r-|\l|$ for $|\la| \leq N+1$ \commentr{ and
$W_{\lambda,\geq k}$ has dimension at most $r-|\lambda|-k$ for
$|\lambda| + k \leq N+1$}
;
\item[(c)] $j$ is sufficiently large (in the sense of
  Lemma~\ref{l:bundle})
for all partitions of integers of length at most $N+1$. 
\end{enumerate}
}

\noindent {\em Proof.}  Let $j_N$ be sufficiently large (in the sense
of Lemma~\ref{l:bundle}) for all partitions of length at most $N+1$
which gives (c).  Part (b) is clear from Lemma~\ref{l:bundle}. 
\commentr{ Note that since $W_{\lambda,\geq k}$
is the image of $W_{\lambda\conc *^k}$ under a projection, we have that
$\dim W_{\lambda,\geq k} \leq \dim W_{\lambda\conc *^k}$.
}
 Taking
$j_N$ even larger, (a) can clearly be satisfied, by ampleness of
$\cL$.  \qed

We will prove Theorem~\ref{T:Bjorn} by showing it modulo ``codimension
$N+1$''  (i.e.\ modulo those classes of dimension at most $-(N+1)$ in
$\cM_{\L}$) for each $N$,  for $j \geq j_N$. 
In what follows, $\cF = \cL^{\otimes j}$ where $j \geq j_N$.

Modulo dimension $< r-N$,
\begin{align*}
W_{\la}
 &\equiv W_{\geq \l} - \sum_{k_1 \geq 1}^{k_1 \leq N-|\la|}  
W_{\l \conc *^{k_1}} 
  \commentr{  - W_{\la  , \geq {N - |\la| + 1} }}\\
 &\equiv W_{\geq \l} - \sum_{k_1 \geq 1}^{k_1 \leq N-|\la|  }
W_{\l \conc *^{k_1}}   \quad \quad \quad   \quad \quad \quad   \quad  \quad \quad   \text{(by Cor.~\ref{c}(b))} \\
& \equiv W_{\geq \l} 
 - \sum_{k_1 \geq 1}^{k_1 \leq N-|\la|}
\commentr{ W_{\geq \l \conc *^{k_1}} }+
\sum_{k_1, k_2 \geq1}^{k_1 + k_2 \leq N - |\la|}  W_{\la \conc *^{k_1}
  \conc \bullet^{k_2}} 
\quad \text{(similarly)} \\
& \equiv W_{\geq \l} - \sum_{k_1 \geq 1}^{k_1 \leq N - |\la| } W_{\geq \la \conc *^{k_1}} +
\sum_{k_1, k_2 \geq1}^{k_1 + k_2 \leq N - |\la| }  W_{\geq \la \conc *^{k_1} \conc \bullet^{k_2}} 
- \sum_{k_1, k_2, k_3 \geq 1}^{k_1 + k_2 +k_3\leq N - |\la|}  W_{\la \conc *^{k_1} \conc
  \bullet^{k_2} \conc 
  \star^{k_3}} 
\\
& \equiv \cdots .
\end{align*}
Continuing (i.e.\ by an easy induction), we may write $W_{\la}$ in
terms of $W_{\geq \la \mu}$ for various $\mu$:
\begin{eqnarray*}
W_{\la}  &\equiv &  \sum_{\mu \in \cQ}^{|\la \mu| \leq N}  (-1)^{ || \mu ||  }W_{\geq \l \mu}
  \\ &\equiv &  \sum_{\mu \in \cQ}^{ |\la \mu| \leq N}  (-1)^{ || \mu ||  } w_{\l \mu} \L^{r
  - |\la \mu| (d+1)} \quad \text{(by \eqref{eq:W})} 
\\ &  \equiv &  \sum_{\mu \in \cQ} (-1)^{ || \mu ||  } w_{\l \mu} \L^{r
  - |\la \mu| (d+1)}  
\pmod {\dim < r-N} %\label{eq:subin}
\end{eqnarray*}
in $\cM$.
We have proved the following.

\tpoint{Proposition}
{\em \label{p:313}For any generalized partition $\lambda$, integer $N\geq|\lambda|$, and $j \geq j_{N}$,\label{p}
$$\frac{ W_{\lambda}} { \L^r} 
\equiv  \sum_{\mu \in \cQ} 
(-1)^{||\mu||} \frac{w_{\lambda\mu}}  { \L^{|\lambda\mu|(d+1)}} \pmod { \codim
  >N}.$$}

Hence modulo codimension $> N$, for $j \geq j_N$, the motivic probability of sections of $\cL^{\otimes j}$ being smooth is
$$
\frac{ W_{\varnothing}} { \L^r} 
=  \sum_{\mu \in \cQ}  (-1)^{ || \mu ||  }  \frac{w_{ \mu}}{ \L^{
    |\mu| (d+1)}} \in  \K.
$$

We have thus shown that for $j \geq j_N$, the left side of
\eqref{eq:Bjorn}  (in the case $s=0$)
stabilizes up to codimension $N$, to the expression of
Proposition~\ref{p} for $\lambda = \varnothing$.
We compare this to the right side of \eqref{eq:Bjorn}.  We have
$$
\frac 1 {\zzeta_X(t)} =  \frac 1 {\sum_{k=0}^{\infty}  \left[ \Sym^k X \right]  t^k}
= \sum_{m=0}^{\infty} \left(  -
\sum_{k=1}^{\infty}  \left[ \Sym^k X  \right] t^k  \right)^{m} = \sum_{\mu \in
\cQ} (-1)^{||\mu||} 
 \left[  \Sym^{\mu} X \right]  t^{|\mu|}.
$$

\point \label{s:s0done}The $s=0$ case of Theorem~\ref{T:Bjorn} is then a consequence of the following proposition.

\tpoint{Proposition}
{\em \label{P:homog}We have
\begin{align*} 
\sum_{\mu \in \cQ} (-1)^{||\mu||} {w_\mu}{t^{|\mu|}}
 &= \sum_{\mu \in \cQ} (-1)^{||\mu||}
\left[ \Sym^{\mu} X \right] t^{|\mu|},
\end{align*}
and hence
\begin{equation}
\frac 1 {\zzeta_X(t)} =  \sum_{\mu \in \cQ} (-1)^{||\mu||}
{w_\mu}{t^{|\mu|}}.
\label{eq:homog}
\end{equation}}

Notice that  $\Sym^{\mu} X$ includes $w_{\mu}$ plus
smaller-dimensional contributions (cf.\ \eqref{e:formal}).  Thus
Proposition~\ref{P:homog}  states that  ``the smaller terms cancel''.

\noindent {\em Proof.}
For $\mu\in\cQ$, we can expand 
$$
\Sym^{\mu} X=\Sym^{f(\mu)} X=\overline{w}_{f(\mu)}=\sum_{\lambda\geq f(\mu)} w_\lambda.
$$
Thus we have
\begin{align*}
\sum_{\mu \in \cQ} (-1)^{||\mu||}
\left[ {\Sym^{\mu} X} \right] {t^{|\mu|}}&= \sum_{\mu \in \cQ} (-1)^{||\mu||} \sum_{\lambda\geq f(\mu)} {w_\lambda}{t^{|\mu|}}  \\
&= 
\sum_{\mu \in \cQ} (-1)^{||\mu||}  {w_\mu}{t^{|\mu|}}
+
\sum_{\mu \in \cQ} (-1)^{||\mu||} \sum_{\lambda> f(\mu)} {w_\lambda}{t^{|\mu|}}.
\end{align*}

We prove the second term in the previous line is $0$ by finding a bijection, for a fixed
$m(\l)$ (hence fixed $w_{\la}$) and fixed $|\mu|$,
between terms with odd $||\mu||$ and terms with even $||\mu||$.  
The map is as follows.  We map a pair $(\lambda,\mu)$ with $\lambda>f(\mu)$ to a pair $(\lambda', \mu')$ with 
$\lambda'>f(\mu')$ that will be constructed below.  

The elements of $\l$ are sums of elements from $f(\mu)$, i.e.\ of $a_1,\dots a_{||\mu||}$.  
Write each element of $\l$ as $a_{i_1} +\dots + a_{i_s}$ with $i_1\leq \cdots \leq i_s$.
We say such an element has length $s$.
Among the longest elements of $\l$, take the lexicographically first sum, call them the \emph{top sums} (there may be a tie) and suppose it ends with $a_k$ (i.e.\ includes at least one $a_k$ term and no $a_j$ term for $j>k$).
If each of the top sums has exactly one $a_k$, and there are no other $a_k$'s in any other terms of $\l$ (which, together, in particular implies that $k\geq 2$), then
we are in {\bf case 1}.  Otherwise, (if a top sum has at least 2 $a_k$ terms, or there are non-top-sum elements of $\l$ containing an $a_k$ summand) we are in {\bf case 2}.

If we are in case 1, note that $k\geq 2$.
As a first attempt, we construct $\lambda'$ by turning all of the $a_k$'s in $\l$ into $a_{k-1}$'s, and we construct $\mu'$ by changing all the $k$'s in $\mu$ into $(k-1)$'s.  However, with this construction $\mu'$ would not necessarily be a partition 
composed of consecutive integers starting with $1$.
So in fact, for each $j\geq k$, we  turn all of the $j$'s appearing in $\mu$ and $\l$ to $(j-1)$'s (as elements or subscripts)
to obtain $\mu'$ and $\l'$, respectively.  We have $\lambda' >f(\mu')$ and $||\mu'||=||\mu||-1$. 
 Since the top sum was lexicographically first among the longest sums, when we replaced $a_k$ by $a_{k-1}$ we do not make any elements of $\l$ equal that were not equal before, and thus $m(\lambda)=m(\lambda')$.  Clearly, $|\mu'|=|\mu|$.
 Furthermore, note that the top sums of $\l$ were changed into elements of $\l'$ which are now top sums of $\l'$.  However, since there are had to be an $a_{k-1}$ appearing somewhere in the original $\l$, the new $(\l',\mu')$ we have created is in case 2.  

In case 2, for all $j>k$ we change all the $j$'s in $\lambda$ and $\mu$ (including in subscripts) to $j+1$'s.  Then we also change one $a_k$ in each top sum  to an $a_{k+1}$ to obtain $\lambda'$, and we change the same number of $k$'s from $\mu$ into $(k+1)$'s to obtain $\mu'$.  
 We have $\lambda' >f(\mu')$ and $||\mu'||=||\mu||+1$. 
Since we change all of the top sums in the same way, we don't make any terms if $\l$ unequal that were previously equal, and thus $m(\lambda)=m(\lambda')$.  Clearly, $|\mu'|=|\mu|$.
  Furthermore, note that the top sums were changed into elements of $\l$ which are now top sums of the new $\l'$.  However, since the $a_{k+1}$'s appear in $\l'$ only in the top sums, and only once in each top sum, the new chain we have created is in case 1. If we then applied the map again, we can see we will get back to $(\lambda,\mu)$.  Similarly, we can check that if we apply the map twice to a $(\lambda,\mu)$ in case 1, we also get back to the original chain.
Thus this operation is an involution, and thus a bijection between even $||\mu||$  and  odd $||\mu||$  terms.
\qed

\bpoint{Proof of Theorem~\ref{T:Bjorn} in general ($s$ singularities)}\label{S:Hyperssing}
In analogy with the inverse of the motivic zeta function
$\zzeta_X^{-1}(t)$, and in light of Proposition~\ref{P:homog}, we
define \highlight{
\begin{equation}\zzeta_{X,\lambda}^{-1} (t) :=\sum_{\mu \in \cQ} 
(-1)^{||\mu||} w_{\lambda \conc \mu} t^{|\lambda \conc \mu|} \in \cM[[t]],\label{eq:b}
\end{equation}
}(a function of $X$).  For example, $\zzeta^{-1}_{X,\varnothing} (t) =
\zzeta_X^{-1} (t)$
by \eqref{eq:homog}.

We will deduce Theorem~\ref{T:Bjorn} from  Proposition~\ref{P:Zform}
and Lemma~\ref{L:norelations}, both  of which will require some time to
prove.   We combine them to prove Theorem~\ref{T:Bjorn} in \S \ref{pf112}.

\tpoint{Proposition}  {\em \label{P:Zform}There exist 
$c_{\pi,\lambda} \in \Z[[t]]$
(independent of $X$ and $\k$), 
%for $\pi \leq m(\la)$, where $b_{\pi,\lambda}$ has constant coefficient
%1,
 such that  }
\begin{equation}\label{E:Zlambda}
\zzeta_{X,\lambda}^{-1}(t)= \sum_{\substack{\pi\in \mathcal{P}\\|\pi|=|\lambda|}}
c_{\pi, \la} w_{\pi} \zzeta^{-1}_X(t). \end{equation}

\noindent {\em Proof.}
Suppose $\la = \prod a_i^{r_i}$, and $\mu$ is a partition whose parts
are distinct from the $a_i$.  Then we have a {\em product rule}
\eqref{eq:productrule} for
$w_{\la \cdot  \mu}$ in terms of $w_{\la} w_{\mu}$ and ``lower order terms''.
Clearly $w_{\la} w_{\mu}$ (interpreted as configurations of distinct
points labeled by $\la$ and distinct points labeled by $\mu$) can be
interpreted as the union (or sum) of $w_{\la \mu}$ along with loci
where some of the points of $\la$ overlap with some of the points of
$\mu$.  We thus have the following formula, where $\mu(i)$ is a
subpartition of $\mu$
corresponding to  which points of $\mu$ overlap with the $a_i$-labeled
points of $\la$. 
\begin{equation}\label{eq:productrule}
\left[ w_{\lambda \conc \mu} \right] = \left[ w_{\lambda} \right]
\left[ w_{\mu} \right] -\sum_{ \substack{
     \mu(i) \text{ not all empty} \\  \cup_i \mu(i) \subset \mu \\
     |\mu(i)|\leq r_i} 
}
\left[ w_{ \left(  a_1^{r_1-|\mu(1)|} a_2^{r_2-|\mu(2) |}\cdots
  \right)  \conc \left(  \mu(1) \mu(2) \cdots \right)\conc \left(
    \mu\setminus \cup \mu(i) \right)} \right] . 
\end{equation}
Now we sum the product rule over all $\mu\in\cQ$ 
to obtain
\begin{eqnarray*}
\zzeta_{X,\lambda}^{-1}(t)&=&\sum_{\mu \in \cQ} (-1)^{||\mu||} \left[ w_{\lambda \conc
  \mu}t^{|\lambda \conc \mu|}  \right] \\
&=& \sum_{\mu \in \cQ} (-1)^{||\mu||} \left[ w_{\lambda} \right] \left[
  w_{\mu}t^{|\lambda\mu|} \right] \\
& & 
-\sum_{\mu \in \cQ} 
\sum_{ \substack{
     \mu(i) \text{ not all empty} \\  \cup_i \mu(i) \subset \mu \\
     |\mu(i)|\leq r_i} }
 (-1)^{||\mu||} 
\left[ w_{ \left(  a_1^{r_1-|\mu(1)|} a_2^{r_2-|\mu(2) |}\cdots
  \right) \conc  \left(  \mu(1) \mu(2) \cdots \right) \conc\left(
    \mu\setminus \cup \mu(i) \right)}
t^{|\lambda\mu|} \right] \\
&=& \left[ w_{\lambda} \right] t^{|\lambda|} \zzeta^{-1}_X(t)
\\ & & -\sum_{\mu \in \cQ} 
\sum_{ \substack{
     \mu(i) \text{ not all empty} \\  \cup_i \mu(i) \subset \mu \\
     |\mu(i)|\leq r_i} } (-1)^{||\mu||}
\left[ w_{ \left(  a_1^{r_1-|\mu(1)|} a_2^{r_2-|\mu(2) |}\cdots
  \right) \conc   \left(  \mu(1) \mu(2) \cdots \right) \conc\left(
    \mu\setminus \cup \mu(i) \right)} \right]
t^{|\lambda\mu|} \quad \text{(by \eqref{eq:homog}).}
\end{eqnarray*}
In this sum, we will group together all the terms where the $\mu(i)$
have some fixed multiplicity sequence $\si(i)$, which is an ordered partition.
Let $F(\si(i))$ be a fixed partition with multiplicity sequence $\si(i)$.
For typographical simplicity, let 
$B = \left( a_1^{r_1-\sum \si(1)} a_2^{r_2-\sum \si(2)}\cdots  \right) \conc
 \left( F(\si(1))  F(\si(2))  \cdots \right)$.

\begin{eqnarray}
\zzeta_{X,\lambda}^{-1}(t)
&=& \left[ w_{\lambda} \right] t^{|\lambda|} \zzeta_X^{-1}(t) \nonumber \\
& & -
 \sum_{\substack{  \si(i) \\ 0\leq' \sum \si(i) \leq r_i}}
\sum_{\mu \in \cQ} \sum_{\substack{  \mu(i) \\ \cup \mu(i) \subset \mu
    \\ m(\mu(i)
  )=\si(i)} }
(-1)^{||\mu||} \left[ w_{\left( a_1^{r_1-|\mu(1)|}
    a_2^{r_2-|\mu(2)|}\cdots \right) \conc 
  \left( \mu(1) \mu(2) \cdots \right)\conc \left( \mu\setminus\cup
    \mu(i) \right)} \right] t^{|\lambda\mu|} \nonumber \\
&=& \left[ w_{\lambda} \right] t^{|\lambda|} \zzeta_X^{-1}(t)  \nonumber \\
& & -
 \sum_{\substack{\si(i) \\ 0\leq'\sum \si(i) \leq r_i}}
t^{\left| \sum_i (\sum \si(i)) \right|}
\sum_{\pi\in\cQ}
(-1)^{||\pi||} \left[ w_{{B\conc\pi}} \right] t^{|B \pi|}
\sum_{\mu \in \cQ} 
\sum_{ \substack{     \mu(i) \\ \cup\mu(i)\subset \mu \\
    m(\mu(i))=\si(i) \\ \mu \setminus\cup \mu(i)\sim\pi}}
(-1)^{||\mu||-||\pi||},\label{eq:new}
\end{eqnarray}
where  $0 \leq' \si(i)$ means that  not all $\sum \sigma(i)$ may be $0$,
and 
 $\sim$ stands for isomorphism of partitions with linearly ordered elements.

We now apply the following lemma, whose proof we defer for a few paragraphs.

\tpoint{Lemma}  {\em \label{L:counting}Given ordered partitions $\si(i) = [\si(i)_j]_j$,  with multiplicity sequence $a(i) = [a(i)_j]_j$, and a partition $\pi$ with linearly ordered elements, we have 
$$\sum_{\mu \in \cQ} 
\sum_{\substack{  \mu(i) \\ \cup\mu(i)\subset \mu \\  m(\mu(i))=\si(i)
    \\  \mu \setminus\cup \mu(i)\sim\pi}}
(-1)^{||\mu||-||\pi||}=\prod_i (-1)^{|\si(i)|}\frac{|\si(i)|!}{\prod_j a(i)_j!}.
$$}

Thus we have
\begin{align}\label{E:rec}
\zzeta_{X,\lambda}^{-1}(t)
&=\left[ w_{\lambda} \right] t^{|\lambda|} \zzeta_X^{-1}
-\sum_{\substack{\si(i) \\ 0\leq'\sum \si(i) \leq r_i} } \left(\prod_i (-1)^{|\si(i)|}\frac{|\si(i)|!}{\prod_j a(i)_j!}\right)
t^{ \left| \sum_i \sum \si(i) \right|}
\zzeta^{-1}_{X,B}(t).
\end{align}
Note that  $m(B)\leq m(\lambda)$ in the merge ordering.
We assume that we know $\zzeta^{-1}_{X,B}$ inductively for $m(B)< m(\lambda)$.
 If $m(B)=m(\lambda)$, then
$\zzeta^{-1}_{X,B}=\zzeta_{X,\lambda}^{-1}$, and we may collect those terms of
the left (note that they all have a nonzero power
of $t$ appearing with them), and solve for $\zzeta_{X,\lambda}^{-1}$.
In particular, we can prove Proposition~\ref{P:Zform} by inducting on
the  level of refinement.

 To conclude the proof of Proposition~\ref{P:Zform}, we prove Lemma~\ref{L:counting}.

\noindent {\em Proof of Lemma~\ref{L:counting}.}
Let $\ell$ be the number of partitions $\si(i)$.
We see that
\begin{eqnarray*}
& & \sum_{\mu \in \cQ} 
\sum_{
\substack{
\mu(i) \\ \cup\mu(i)\subset \mu \\ m(\mu(i))=\si(i) \\ \mu \setminus\cup \mu(i)\sim\pi}}
(-1)^{||\mu||-||\pi||}\\
& =&\sum_{ \mu , \mu_1 \dots,\mu_{\ell-1} \in \cQ} 
\sum_{  \substack{  \mu(\ell) \subset \mu_{\ell-1} \\
    m(\mu(\ell))=\si(\ell) \\
\mu_{\ell-1}\setminus\mu(\ell)\sim \pi }
}
(-1)^{||\mu_{\ell-1}||-||\pi||}\cdots \sum_{ \substack{  \mu(2)\subset \mu_1 \\ m(\mu(2))=\si(2) \\
\mu_1\setminus\mu(2)\sim \mu_2
}}(-1)^{||\mu_1||-||\mu_2||} 
\sum_{\substack{  \mu(1)\subset \mu \\ m(\mu(1))=\si(1) \\
    \mu\setminus\mu(1)\sim \mu_1 } }
(-1)^{||\mu||-||\mu_1||} 
\end{eqnarray*}
and so we can reduce to the case in which $\ell=1$, i.e.\ there is only one partition $\si(i)$, which we call $\sigma$, with
multiplicity sequence $[a_j]_j$.

Fix an integer $k$ and consider the case when $||\mu||-||\pi||=k$.
We will start with $\pi$ and need to count how many $\mu$ with
$m(\mu)=\sigma$ we can add to $\pi$ to obtain a partition with $||\pi||+k$ linearly ordered distinct elements.
To choose where the new $k$ elements go in the ordering, there are
$\binom{||\pi||+k}{k}$ possible choices.  Then there are 
$\binom{||\pi||}{|\sigma|-k}$ choices for which elements of $\pi$ will also be elements of $\mu$.
Once we have made those choices, there are $\frac{|\sigma|!}{\prod_j a_j!}$ choices for
how to assign the $|\sigma|$ multiplicities in $\sigma$ to these $|\sigma|$ locations.
The well-known identity 
$$
\sum_{k=0}^{|\sigma|} (-1)^k \binom{||\pi||+k}{k}\binom{||\pi||}{|\sigma|-k}=(-1)^{|\sigma|}.
$$
(which can be proved with generating functions, for example), concludes the proof
 of Lemma~\ref{L:counting}.
 \qed

This in turn concludes the proof of 
 Proposition~\ref{P:Zform}.
\qed
 
\tpoint{Lemma (``no unexpected universal linear relations among the
  $w_{\lambda}$'')}
{\em \label{L:norelations}The relations $w_\lambda=w_{\lambda'}$ for
$m(\lambda)=m(\lambda')$ generate all the $\Z$-linear relations among
the $w_\lambda$ that hold for {\em all} smooth
projective varieties 
  $X$ of pure dimension  $d>0$ over  {\em all} fields $\k$.}

We conclude the proof of
Lemma~\ref{L:norelations}
in \S \ref{pf:l38}.

\tpoint{Lemma (``no unexpected universal relations among symmetric powers'')}
{\em \label{L:SymInd}Suppose we have $f \in \Z[x_1, x_2, \dots]$ such that $f([\Sym^1 X],
  [\Sym^2 X], \dots) = 0 \in \Kvar$ for all smooth projective varieties $X$ of pure dimension over all
  fields $\k$.  Then $f \equiv 0$.}

In fact, the argument uses only dimension $1$.

\noindent {\em Proof (Poonen, cf.\ \cite[\S 3.3]{Poonen04}).}
Suppose we have a non-zero $f \in \Z[x_1,
\dots, x_n]$ (where $x_n$ appears in $f$) such that 
$f(\Sym^1X, \dots, \Sym^nX)=0$ for all varieties $X$ over any
field $\k$.

Suppose $\k = \FF_q$.  The information of $\#(\Sym^1
X)$, \dots, $\#(\Sym^n X)$ is equivalent to the information of the
number of points of degree $1, \dots, n$ on $X$.  Thus there exist
non-negative integers $a_1, a_2, \dots, a_n$ such that there does not
exist an $X$ with $a_i$ points of degree $i$.
%  Then in a large enough extension $\k$ of $\FF_p$, 
%requiring that $X$ have $a_i$ points of degree $i$ for $1 \leq i \leq n-1$ imposes a nontrivial constraint
%on the number of points of degree $n$ (i.e.\ there are only a finite number of possibilities).
%(Translation:  we've got a nontrivial polynomial in n variables in which the nth variable actually appears.   If we specialize the first n-1 %variables in a general enough way, then we still have a nonzero polynomial in the nth variable.)
%Now let $a_i$ be positive integers such that $f(a_1, \dots,  a_n)\neq 0$.

We construct a smooth projective curve with $a_i$ points of degree $i$
for $i=1, \dots, n$, yielding a contradiction.  Choose $N$ large
enough so that $\PP^N_{\FF_q}$ has at least $a_i$ points of degree
$i$, and pick $a_i$ points of degree $i$ in $\PP^N_{\FF_q}$.  Then a
standard argument shows that we can find $N-1$ hypersurfaces $f_1$,
\dots, $f_{N-1}$ intersecting completely (i.e.\ intersecting in a
curve), such that of the points of degree at most $n$, the $f_i$ pass
precisely through our chosen points; and such that the $f_i$ are
linearly independent in the tangent space at each of our chosen
points.  The complete intersection of the $f_i$ is a curve containing
precisely the desired number of points of degree at most $n$, and
smooth at those points.  Take the normalization of this curve (which
will not introduce any more points of small degree).  \qed
 
\epoint{Observation} \label{o:chains} We have the following universal formula for
$w_{\l}(X)$ in terms of the symmetric powers of $X$.
\begin{equation}\label{eq:chains}
\left[ w_{\l}(X) \right] = \sum_k\sum_{\l=\mu_0\ll \mu_1 \ll \dots \ll\mu_k} (-1)^k
\left[ \Sym^{m(\mu_k)} X \right] .
\end{equation}
(By Lemma~\ref{L:SymInd}, this formula is unique.)
Here (for the purpose of this argument only)
we say $\lambda \ll  \lambda'$ if $f(\lambda)<\lambda'$.
We show \eqref{eq:chains} by induction on the length of $\lambda$.  Clearly it is true for $|\lambda|=1$.
Then,
\begin{eqnarray*}
\left[ w_{\l} (X) \right] &=& \left[  w_{f(\l)} (X)  \right] \\ &= &
\left[ \wbar_{f(\l)}(X) \right] - \sum_{f(\l)<\mu} \left[ w_\mu
(X) \right]  \\  &=& \left[ \Sym^{m(f(\l))}  X \right] -  \sum_{f(\l)<\mu}  \sum_{\mu=\mu_0\ll
  \mu_1 \ll  \dots \ll \mu_k} (-1)^k \left[ \Sym^{m(\mu_k)} X \right] .
%&=\Sym^{m(\l)} -  \sum_{\l\ll p}  \sum_{[p=p_0\ll  p_1 \ll  \dots \ll p_k]} (-1)^k \Sym^{m(p_k)}\\
%&=\Sym^{m(\l)} +  \sum_{\l\ll p_0\ll  p_1 \ll  \dots \ll p_k, k\geq 0} (-1)^{k+1} \Sym^{m(p_k)}\\
\hidden{&=\sum_{\l=\mu_0\ll  \mu_1 \ll  \dots \ll \mu_k} (-1)^k \Sym^{m(\mu_k)}}
\end{eqnarray*}

\epoint{Proof of Lemma~\ref{L:norelations}}
\label{pf:l38}If we had a finite non-trivial relation $\sum c_\lambda w_\lambda=0$
only involving terms with distinct multiplicity sequences, then we can use \eqref{eq:chains} to write it
as an algebraic relation on $[\Sym^i X]$, holding for all smooth projective $X$ of pure dimension over
all fields $\k$.  The terms with non-zero $c_\lambda$ for maximal $|\lambda|,$
will give terms $c_\lambda [\Sym^{\lambda} X]$, which will be the maximal degree terms in the algebraic relation
(where $[\Sym^i X]$ has degree $i$).  Thus we will obtain a non-trivial algebraic relation
$f([\Sym^1 X], ..., [\Sym^j X] )=0$ for some $j$, contradicting Lemma~\ref{L:SymInd}.
\qed

\bpoint{Proof of Theorem~\ref{T:Bjorn}}
\label{pf112}We finally prove
Theorem~\ref{T:Bjorn} using Proposition~\ref{P:Zform} and Lemma~\ref{L:norelations}.

From \eqref{eq:b},  $\sum_{s\geq 0} \zzeta_{X,*^s}^{-1}(t) =1$:  consider
the contribution from the right side of \eqref{eq:b} to the term $w_{\nu}$ for each
$\nu$, and note that if $\nu \neq \varnothing$, 
say $\nu = *^k \pi$, then  $w_\nu$ has contributions from the $s=0$
term and the $s=k$ term, with opposite signs.

 Thus, multiplying \eqref{E:Zlambda} by $\zzeta_X(t)$, and summing over $\lambda=*^s$ for $s\geq 0$, we obtain
$$
\sum_{\lambda\in \mathcal{P}} w_\lambda t^{\Sigma \la}=\zzeta_X(t)
=\sum_{s\geq 0} \zzeta_{X,*^s}^{-1}(t) \zzeta_X(t)=\sum_{\pi\in\mathcal{P}} {c_{\pi,*^{|\pi|}}}
w_{\pi}.
$$
Since Lemma~\ref{L:norelations} implies that the linear relations among the $w_\lambda$ with fixed
$|\lambda|$ generate all the linear relations among the $w_\lambda$, we can deduce that
\begin{equation}\label{E:Zs}
 \zzeta_X^s(t)=\sum_{\substack{\lambda\in \mathcal{P}\\|\lambda|=s}} w_\lambda t^{\Sigma \la}
=\sum_{\substack{\pi\in\mathcal{P}\\|\pi|=s} }{c_{\pi,*^{|\pi|}}}w_{\pi}= \zzeta_{X,*^s}^{-1}(t) \zzeta_X(t).
\end{equation}
Thus for $j \geq j_s$,
\begin{eqnarray*}
\frac{ W_{*^s}} { \L^r} & \equiv &\
\sum_{\mu \in \cQ} (-1)^{||\mu||}  \frac {w_{\mu *^s}} { \L^{|\mu
    *^s|(d+1)}} \pmod {\codim >s}
\quad \text{(Prop.\ \ref{p:313})} \\ 
&=&  \zzeta_{X,*^s}^{-1} \left(  \L^{-(d+1)} \right) \quad \quad \quad \quad
\quad  \quad \quad \quad  \quad \quad \quad \quad \quad \quad 
\text{(from \eqref{eq:b})} \\
&=& 
\left( \sum_{|\lambda|=s} \frac{ w_\lambda }  { \L^{(d+1)\sum
      \lambda} } \right) \frac 1 {\zzeta_X(\L^{-(d+1)})}
\quad  \quad \quad \quad \quad \quad \quad  \text{(from  \eqref{E:Zs}),}
\end{eqnarray*}
which proves Theorem~\ref{T:Bjorn}.

\epoint{Remark:  $s$ ordered points (cf.\ Remark~\ref{r:Bjorn}(vii))}
\label{r:sordered}Let $\la$ be the partition $1
2 \cdots s$.  Then 
$$\frac {  [H^0(X,\cL^{\otimes j})^\text{$s$ ordered}]}  {[
  H^0(X,\cL^{\otimes j})]} = \frac {[W_{\la} ]} { \L^r}$$ by the definition of
$W_{\la}$ (where $r = h^0(X,
\cL^{\otimes j})$ as in Lemma~\ref{l:bundle}).  Then by
Proposition~\ref{p:313},
$$\lim_{j \rightarrow \infty} \frac {  [H^0(X,\cL^{\otimes j})^\text{$s$ ordered}]}  {[
  H^0(X,\cL^{\otimes j})]} =   \sum_{\mu \in \cQ} 
(-1)^{||\mu||} \frac{[w_{\lambda\mu}]}  { \L^{|\lambda\mu|(d+1)}},$$
which  is $\zzeta_{X, \la}^{-1}(\L^{-d-1})$ by the definition \eqref{eq:b}
of $\zzeta_{X, \la}^{-1}(t)$.  By inserting Lemma~\ref{L:counting} into
\eqref{E:rec},  after slight rearranging we have 
$$\zzeta_{X, \lambda}^{-1}(t) = 
[w_{\lambda}]t^{|\lambda|} \zzeta_X^{-1}(t)
-\sum_{k=1}^{s} \binom{s}{k}(-1)^{k} t^{k}
\zzeta_{X,\lambda}^{-1}(t)$$
from which $\zzeta_{X,\lambda}^{-1}(t)=\frac{w_{\lambda}t^{s}
  \zzeta^{-1}_X(t)}{(1-t)^{s}}$. Remark~\ref{r:Bjorn}(vii)  follows.

\section{Motivic stabilization of symmetric powers}\label{S:MSSP}

In this section, we prove three statements given  in the introduction.

\tpoint{Proposition} {\em 
Suppose $X$ is a geometrically irreducible variety.   Then $[\Sym^n (X \times
\A^1)] = \L^n \times [\Sym^n X]$, and hence   motivic stabilization
of symmetric powers holds for $X$ if and only if it holds for $X \times \A^1$.
\label{p:XL}\label{p:A}}

This is essentially  \cite[Statement 3]{GZLMH04}, and also follows by applying Totaro's argument of 
\cite[Lemma 4.4]{Got01}.   The same argument applies for MSSP${}_\phi$
for any motivic measure $\phi$.   \cut{
\begin{theorem}[{\bf Totaro's Lemma~\cite[Lemma 4.4]{Got01}}]
$\Sym^n \A^m$ can be cut-and-pasted to $\A^{nm}$, so  
$\Sym^n \A^m = \A^{nm}$ in $\cM$.\label{t:totaroslemma}
\end{theorem}

  (We find it remarkable that this fact was discovered so late.   Totaro's proof is very short but nonconstructive.   An explicit cutting and pasting is implicitly given in \cite[Thm.~1]{GZLMH11}, see also \cite[Statement~2]{GZLMH04}.)}

\tpoint{Proposition} {\em 
Suppose $X$ is a geometrically irreducible variety, $U \subset X$ is a dense open set, and $Y \subset X$ is the complementary closed set.    Then 
 \begin{equation}\label{eq:DD}\lim_{n \rightarrow \infty} \frac{ \left[ \Sym^n X \right] }  {\L^{dn}} = \zzeta_Y( \L^{- d})
 \lim_{n \rightarrow \infty} \frac{ \left[ \Sym^n U \right] }{ \L^{dn }}. 
\end{equation}
 More precisely the limit on the left exists (motivic stabilization
 of symmetric powers holds for $X$)  if
 and only if the limit on the right exists (motivic stabilization of
 symmetric powers holds for $U$), and in this case \eqref{eq:DD} holds.
 In particular, if $X_1$ and $X_2$ are birational geometrically irreducible
 varieties, then motivic stabilization of symmetric powers holds for $X_1$ if and only if it holds for $X_2$.\label{p:motivicbirational}\label{p:B}}
 
 \noindent {\em Proof.}
 We prove the result modulo dimension $-k$ classes, by induction on $k$.   The case $k=0$ is trivial.
For all $n$,
 \begin{equation}\label{eq:stab}\frac{ \left[ \Sym^n X  \right]}
   {\L^{dn} }= \frac {\left[ \Sym^n U \right] } {\L^{dn}} +
 \frac { \left[ \Sym^{n-1} U  \right]} {\L^{d(n-1)}}  \frac { \left[
     \Sym^1 Y \right]} {
   \L^{d}} + \cdots  +
 \frac { \left[ \Sym^{n-k} U  \right] } {\L^{d (n-k) }}  \frac {
   \left[ \Sym^k Y \right]} { \L^{d k}}
 \end{equation}
 modulo classes of dimension less than $-k$.  (Here we use that $\dim
 [ \Sym^m Y ] / \L^{dm } \leq -m$, and 
 $\dim [ \Sym^m ] U / \L^{d m} = 0$.)    If the symmetric powers of $U$
 motivically stabilize, to $\Sinf(U)$ say, then the right side
 stabilizes to $\Sinf(U) \zzeta_Y(\L^{- d})$ modulo classes of
 dimension $< -k$, as desired.  On the other hand, if  the symmetric
 powers of $X$ stabilize, and the symmetric powers of $U$ stabilize up
 to dimension $-k+1$, then everything  in \eqref{eq:stab} stabilizes
  (modulo classes of dimension $<-k$) except for possibly $
 {[ \Sym^n U]} /  { \L^{dn}}$; but then this class must stabilize as well, as desired.
 \qed 

\tpoint{Proposition} {\em
If $X$ is a geometrically irreducible smooth projective curve with a
rational point, then motivic
stabilization of symmetric powers holds for $X$.\label{p:curve}\label{p:C}}

\noindent {\em Proof.} For $n > 2g-2$ (where $g$ is the genus of $X$),
$\Sym^n X$ is a (Zariski)  $\PP^{n-g}$-bundle over $\operatorname{Jac}
C$ (where we use the point to determine an isomorphism
$\operatorname{Pic}^n X \cong \operatorname{Jac} C$, and so that
$\Sym^n X$ is a Zariski bundle), so  $[\Sym^n X]
= [\PP^{n-g}] [ \operatorname{Jac} C]$. \qed

\section{Moduli of points:  Configuration spaces}

\label{s:cs}Throughout this section, $X$ will be a geometrically irreducible
variety of dimension $d$.  For notational convenience, define
\highlight{$\M = \L^d$}.  In this section, we will show that (under
appropriate motivic stabilization hypotheses) the classes of all
``discriminants'' ``stabilize'' (as the number of points tends to
$\infty$) to finite formulas in terms of motivic zeta values, which can be interpreted in terms of ``motivic probabilities''
(Theorem~\ref{T:IntroPointsLimits}, see Corollary~\ref{C:stabform}
below).   Corollary~\ref{C:stabform} is a consequence
of an unconditional statement about generating functions
(Theorem~\ref{T:Kform}), which is the most complicated result in this
section.  Explicit special cases can be shown without the full
strength of Theorem~\ref{T:Kform}, and are sprinkled
throughout.

%\bpoint{Key results}\label{S:genfunc}
We first  name the generating functions that are the subject of our investigation.
For each partition of  positive integers $\nu$, define 
\highlight{$$
K_{1^{\bullet}\nu}(t):=\sum_{j} w_{1^j \nu}(X) t^j \quad \textrm{and} \quad
\Kbar_{1^\bullet \nu}(t):=\sum_{j} \wbar_{1^j \nu} (X) t^j.
$$}
In the proofs of our results, we use a generalization of
$K_{1^{\bullet}\nu}$(t).  For an integer $a$, define
\highlight{$w_{(<a \vdash j)\nu}
  (X):=\sum w_{\mu\nu} (X)$}, where the sum is over the set of
partitions $\mu$ of $j$ with all parts less than $a$.
 Define
\highlight{\begin{equation}\label{eq:defK}
K_{(<a)\nu}:=\sum_{j} w_{(<a \vdash j) \nu} (X) t^j,
\end{equation}}
so $K_{(<2)\nu}=K_{1^{\bullet}\nu}$.
Informally: 
$(<a \vdash j)$ refers
the set of partitions of $j$ with parts $<a$, and 
 $(<a)$ refers to the set of all partitions with all parts
$<a$.

\tpoint{Theorem}\label{T:Kform} 
{\em 
For each partition $\nu$ of positive integers and each $a\geq 2$,
there exist  universal formulas $A_{\nu,a}(t)$, $B_{\nu,a}(t)$, $C_\nu(t)$, and
$D_\nu(t)$ (recursively defined in Propositions~\ref{P:Krecursion} and~\ref{P:K}(b)) such that
\begin{itemize}
 \item  $A_{\nu,a}(t)$ is a $\Z[t]$-linear combination of 
   $w_{\nu'}$, where $|\nu'|=|\nu|$ and $m(\nu')\leq m(\nu)$;
\item $C_\nu(t)$ is a $\Z[t]$-linear combination of 
$w_{\nu'}  / \zzeta_X(t^{i})$, where 
 $|\nu'|\leq |\nu|-1$, $i\in \nu$, and $i\geq 2$;
\item  and $B_{\nu,a}(t), D_\nu(t) \in \Z[t]$, both having constant coefficient $1$;
\end{itemize}
such that for any $X$  and  $\k$, 
\begin{enumerate}
\item[(a)] 
when $\nu$ has all
parts at least $a$,
$$K_{(<a)\nu} (t) =\frac{\zzeta_X(t) }{\zzeta_X(t^{a})
}\frac{A_{\nu,a}(t)}{B_{\nu,a}(t)},$$
\item[(b)]
$$\Kbar_{1^{\bullet} \nu} (t) = \zzeta_X(t) \frac{C_\nu(t)
}{D_\nu(t)}.$$
\end{enumerate}
}

These formulas are independent of $X$ and $\k$, in the sense that
their only dependence on $X$ is via the universal formulas 
(in
terms of the $\Sym^{\bullet} X$) for 
$w_{\nu'}(X)$ (given in \S \ref{o:chains}) and $\zzeta_X(t)$.  
Because the formulas for $A_{\nu, a}(t)$ through $D_{\nu}(t)$ are
finite, the formulas for $K_{(<a) \nu}(t)$ and $\Kbar_{1^{\bullet} \nu}(t)$ have
finite descriptions in terms of zeta functions.  
Part (a)  will follow from
Proposition~\ref{P:Krecursion}, and part (b) will follow from 
part (a) and~\ref{P:K}(b) (see the comment
after the statement of Proposition~\ref{P:K}).    Before embarking on
the proof, we give
some interpretations and special cases.

\bpoint{Interpretations and consequences}

By the following lemma, the  formulas for  $K_{(<a)\nu}(t)$ (hence 
$K_{1^{\bullet} \nu}(t)$), and $\Kbar_{1^{\bullet} \nu}(t)$
 imply similar formulas
for the limits of their coefficients.

\tpoint{Lemma} \label{L:limits}{\em 
Suppose 
$Y(t) = E(t) \zzeta_X(t)$,
where $Y(t)  =\sum_{j} Y_j t^j$ and $E(t) = \sum_j E_j t^j$ both lie in
$\mathcal{M}_\L[[t]]$.
\begin{enumerate} \item[(a)]
If $X$ satisfies  MSSP${}_\phi$, and
$E(\M^{-1})$ exists (i.e., converges in $\phi(\K)$), then
\begin{equation}
\label{eq:july1}\lim_{j\ra\infty} \frac{Y_j}{\M^j} = E(\M^{-1}) \Sinf_\phi (X) \quad
\text{in  $\phi(\K)$.}
\end{equation}
If further the $[\Sym^j X]$ are invertible in $\phi(\K)$ (e.g.\ if $X$
is rational or $\phi = HS$, \S \ref{o:Syminv}), then
\begin{equation}\label{eq:july2}\lim_{j\ra\infty} \frac{Y_j}{ \left[
      \Sym^j X \right]} = E(\M^{-1}) \quad
\text{in  $\phi(\K)$.}\end{equation}
\item[(b)]
 If $\k = \F_q$, and $\# E(  1/ q^d)$ converges (in
  $\R$), then 
$$\lim_{j\ra\infty} \frac{\# Y_j}{ \# \Sym^j X} = \# E(1/q^d).$$\end{enumerate}
}

The stabilization of the symmetric powers of
$X$
in the motivic measure $\phi$, that is,
$\Sinf_\phi(X) = 
\lim_{n\ra\infty} \frac{\left[ \Sym^n X \right]}{\L^{dn}}$,  was defined in \S \ref{s:mssp1}.
Part (b) may be interpreted as  inspiration for (a), but 
does not follow from (a), as the point-counting map $\#$ does not extend to a map $\K
\rightarrow \R$ (see \S \ref{s:dimfil}).

\noindent {\em Proof.} 
(a) We prove \eqref{eq:july1} in the case where $\phi$ is the identity, i.e.\ when
motivic stabilization holds in general.  The extension to a particular
$\phi$ is straightforward using the continuity of $\phi$, but the
argument is notationally more cumbersome, hence left to the reader.

To prove the result, we prove the result modulo dimension at most
$-b$, for any $b$.  By hypothesis $E_i\M^{-i}$ is bounded above in
dimension by some integer $b_1$.  Let $b_2$ be such that if $i\geq
b_2$, then $E_i\M^{-i}$ has dimension at most $-b$, and $b_3$ be such
that for $i\geq b_3$, we have that $\Sinf(X)- [\Sym^i X] \M^{-i}$ has
dimension at most $-b-b_1$.  For $j\geq b_2+b_3$ we have, modulo
dimension at most $-b$,
$$
 \frac{Y_j}{\M^j} =\sum_{i=0}^j \frac{E_{j-i}}{\M^{j-i}}\frac{
   \left[ \Sym^i X \right]}{\M^i}
\equiv \sum_{i=j-b_2}^j \frac{E_{j-i}}{\M^{j-i}}\frac{ \left[ \Sym^i X \right]}{\M^i}
\equiv \sum_{i=j-b_2}^j \frac{E_{j-i}}{\M^{j-i}}\Sinf(X)
\equiv \sum_{i=0}^\infty \frac{E_{i}}{\M^{i}}\Sinf(X).
$$

Part (b) is an exercise in convergent power series, using the fact (from
the Weil conjectures)  that the unique root of the denominator of
$\zeta_X(t)$ with the smallest absolute value is $1/q^d$ (cf.\
Motivation~\ref{motivation}(i)).
 \qed

 Combining Theorem~\ref{T:Kform} (still to be proved) and Lemma~\ref{L:limits}, we have
 the limits of normalized configuration spaces promised in the
 introduction (Theorem~\ref{T:IntroPointsLimits}). The hypotheses of Lemma~\ref{L:limits} are straightforward to check.

\tpoint{Corollary} {\em \label{C:stabform} 
Suppose $\nu$ is a partition of positive integers, all at least $2$.
\begin{enumerate}
\item[(a)] If $X$ satisfies
  MSSP${}_\phi$, then 
$$
\lim_{j\ra \infty} \frac{w_{1^j\nu}}{\M^{j}} = \frac 1
{\bzeta_X(2d)}    \frac{A_{\nu,2} (\M^{-1})}{B_\nu (\M^{-1})}
\Sinf_\phi (X)\quad \textrm{and} \quad
\lim_{j\ra \infty} \frac{\wbar_{1^j\nu}}{\M^{j}} =
 \frac{C_{\nu,2} (\M^{-1})}{D_\nu (\M^{-1})}  \Sinf_\phi (X)
$$
in $\phi(\K)$.
If furthermore the $[\Sym^j X]$ are invertible (e.g.\ if $X$ is rational
or $\phi=HS$, \S \ref{o:Syminv}), then
$$
\lim_{j\ra \infty} \frac{\left[ w_{1^j\nu} \right] }{ \left[\Sym^{j +
      |\nu|} X \right]} = \frac 1
{\M^{|\nu|} \bzeta_X(2d)}   \frac{A_{\nu,2} (\M^{-1})}{B_\nu (\M^{-1})}
  \quad \textrm{and} \quad
\lim_{j\ra \infty} \frac{ \left[\wbar_{1^j\nu} \right]} { \left[
    \Sym^{j+|\nu| }X \right]} =
\frac 1
{\M^{|\nu|} }  \frac{C_{\nu,2}(\M^{-1})}{D_\nu(\M^{-1})}.
$$
\item[(b)]  If $\k = \F_q$, then 
$$
\lim_{j\ra \infty} \frac{\# w_{1^j\nu}}{\# \Sym^{j + |\nu|} X} = \frac 1
{q^{d|\nu|} \zeta_X(2d)}    \frac{\# A_{\nu,2}
      (q^{-d})}{\# B_\nu (q^{-d})}
  \quad \textrm{and} \quad
\lim_{j\ra \infty} \frac{\# \wbar_{1^j\nu}}{\# \Sym^{j+|\nu| }X} =
\frac 1
{ q^{d|\nu|} }  \frac{\# C_{\nu,2} (q^{-d})}{\# D_\nu (q^{-d})}.
$$
\end{enumerate}}

\bpoint{Determining the universal formulas
$A_{\nu, a}(t)$ through $D_{\nu}(t)$}

We now begin the proof of Theorem~\ref{T:Kform}.  En route,
we in effect prove special cases, giving explicit descriptions of 
$A_{\nu, a}(t)$ through $D_{\nu}(t)$.
We start with an important base case.

\tpoint{Proposition (see \S \ref{promise56})}\label{P:Kbase} {\em If $a>1$, \begin{enumerate}  
\item[(a)] $K_{(<a)}(t)=\zzeta_X(t)/\zzeta_X(t^a)$.
\item[(b)]
 $\Kbar_{1^{\bullet} (a)}(t)=t^{-a} \zzeta_X(t) (1-1/\zzeta_X(t^a))$.\end{enumerate}}

\noindent {\em Proof.}
We have $$\left[ \Sym^j X \right] t^j=
\sum_{\la \vdash j}  \left[ w_{\lambda} \right]
t^{\sum \lambda}
=
\sum_{\mu\in \mathcal P}  \left[ w_{a \times\mu} \right] t^{a\sum \mu}
\left[ w_{(<a \vdash j-a\sum \mu )} \right]t^{j-a\sum \mu},$$
where the notation  $a \times p$ (used only in this proof)  denotes the partition obtained by multiplying all of the elements of $p$ by $a$.
(Given $\la$, to find the $\mu$ on the right side, ``round down'' the
parts of $\la$ to the next multiple of $a$, then divide the result by
$a$.  A similar idea, with $a=2$, is used in the proof of
Theorem~\ref{T:Kssing} below.)
Using $[w_{a \times \mu}]=[ w_\mu]$, we have 
\begin{eqnarray*}\zzeta_X(t) &=&
\sum_j \left[ \Sym^j X \right] t^j \\ &=&
\sum_j \sum_{\mu\in \mathcal P} \left[ w_{a \times\mu} \right]
t^{a\sum \mu} \left[ w_{(<a
  \vdash j-a\sum \mu )}  \right] t^{j-a\sum \mu} \\ &=&
 \left( \sum_{\mu \in \mathcal P} \left[ w_\mu \right]  t^{a \sum \mu} \right)
\left(  \sum_k  \left[  w_{(<a \vdash k)}  \right] t^k \right) \\ &= &
\zzeta_x(t^a) K_{(<a)}(t),\end{eqnarray*}
yielding part (a).  
Part (b) then follows from (a) via the identity 
 $K_{(<a)}(t) + t^a \Kbar_{1^{\bullet} (a)}(t)=\zzeta_X(t)$:  each partition of
$n$
either has all parts smaller than $a$, or else at least one part at
least $a$. \qed

\tpoint{Theorem}\label{T:Kssing} {\em 
$$
\sum_j \left[ \Sym^j_s X \right] t^j = \frac{\zzeta^{[s]}_X(t^2)\zzeta_X(t)}{\zzeta_X(t^2)}.
$$
}

Theorem~\ref{T:IntroPointsssing} follows by combining Theorem~\ref{T:Kssing} and Lemma~\ref{L:limits}.

\noindent {\em Proof.}
Let \highlight{$\mathcal{S}_s$} be the set of partitions of positive
integers $\lambda$ with $|\lambda|=s$ and with all parts even.
For a partition $\lambda$ with all parts even, let
\highlight{$\mathcal T_{\lambda,j}$}  be the set of all partitions
$\mu \vdash j$ such that 
 $\{2\lfloor \mu_i/2 \rfloor \; | \;  2\lfloor \mu_i/2 \rfloor>0 \}=\lambda,$
(i.e.\ $\lambda$ is obtained from $\mu$ by rounding the parts down to the nearest even integer and discarding $0$'s, cf.\ the
proof of   Proposition~\ref{P:Kbase}(a)).
(The notation $\mathcal S_s$ and $\mathcal T_{\la, j}$ will only be
used in this proof.)

For a partition $\la$ with all parts even, 
$w_\lambda w_{1^{j-\sum \lambda}}= \sum_{\mu\in \mathcal T_{\lambda,j
  }}w_\mu$, so
$$
\sum_{\lambda\in\mathcal  S_s}  \left[ w_\lambda \right] \left[w_{1^{j-\sum
      \lambda}} \right]=
\sum_{\lambda\in\mathcal  S_s} \sum_{\mu\in\mathcal  T_{\lambda,j }}
\left[ w_\mu \right]
=  \left[ \Sym^j_s X \right] ,
$$
as the middle double sum enumerates the partitions $\mu$ of $j$ with precisely
$s$ multiple points.
Thus 
\begin{align*}
 \sum_j \left[ \Sym^j_s X \right] t^j = \sum_j
 \sum_{\lambda\in\mathcal  S_s}  \left[ w_\lambda \right] \left[
   w_{1^{j-\sum \lambda}} \right] t^j=
 \left( \sum_{\lambda\in\mathcal  S_s} \left[ w_\lambda \right] t^{\sum \lambda}
 \right) \left( \sum_k \left[ w_{1^k} \right]  t^k \right)=\zzeta^{[s]}_X(t^2) \frac{\zzeta_X(t)}{\zzeta_X(t^2)}.
\end{align*}
where the last equality uses
Proposition~\ref{P:Kbase}(a).  \qed

Temporarily  (for the purpose of
Proposition~\ref{P:Krecursion}) define 
\highlight{$\mathcal A_{<a}(\nu)$} as the set of all partitions
obtained by adding an element of $\{ 0, \dots, a-1 \}$  to each of
the parts of $\nu$.  (Think:  ``$\mathcal A$dd $<a$ to the parts of
$\nu$''.)  For
example, $[x+2,x+2,x,y+1,y]\in\mathcal A_{<a}([x,x,x,y,y])$ for
$a\geq 3$. 

\tpoint{Proposition (recursion for $K_{(<a)\nu}(t)$)}\label{P:Krecursion} {\em 
For any formalization $\nu$ of a partition,
$$
 K_{(<a)\nu}(t)= 
\frac {\frac{\zzeta_X(t)}{\zzeta_X(t^a)}{w_\nu}- \sum_{ \substack{ 
\nu' \in \mathcal A_{<a}(\nu)
\\ m(\nu')<m(\nu)  }
}  {K_{(<a)\nu'}(t)}{t^{\sum \nu'-\sum \nu} }}
{
\sum _{ \substack{  \nu' \in \mathcal A_{<a}(\nu)
\\ m(\nu')=m(\nu) 
}}  t^{\sum \nu'-\sum \nu} 
}.
$$
}

Note that the denominator is a polynomial with constant coefficient
$1$, as  $\nu' = \nu$ appears in the bottom sum. Theorem~\ref{T:Kform}(a)
follows inductively from Proposition~\ref{P:Krecursion}, using the base case $\nu=\varnothing$, because we may replace
$\nu$ by its formalization.    (The assumption that all the
parts of $\nu$ are at least $a$ in Theorem~\ref{T:Kform}(a) arises  because of this need to
replace $\nu$ by its formalization.)

\noindent {\em Proof.}
For any $a$ and $j$,
\begin{equation}
\label{eq:oldlem}
w_{(<a \vdash j)}w_\nu  = \sum _{\nu' \in\mathcal A_{<a}(\nu)} w_{{(<a\vdash j-\sum \nu' +\sum \nu)}\nu'}.
\end{equation}
Reason:  when multiplying $w_{(<a \vdash j)}$ with $w_{\nu}$, the
right side keeps track of ``how the points parametrized by $w_{(<a \vdash j)}$ and
$w_{\nu}$ overlap''.    
Multiplying \eqref{eq:oldlem} by $t^j$ and summing over all $j$, we
have
\begin{eqnarray*}
\sum_j \left[ w_{(<a \vdash j)} \right] \left[ w_{\nu}  \right] t^j  &=& 
\sum_{j } t^j\sum _{\nu' \in
  \mathcal A_{<a}(\nu)} \left[ w_{{(<a \vdash j-\sum \nu' +\sum
      \nu)}\nu'} \right]\\
K_{(<a)}(t) \left[ w_{\nu} \right] &=& 
\sum _{\nu' \in
  \mathcal A_{<a}(\nu)}  t^{\sum \nu' - \sum \nu}
\sum_{k } t^k
\left[ w_{(<a \vdash k)\nu'} \right]  \quad \text{(by \eqref{eq:defK})} \\
\frac{\zzeta_X(t)}{\zzeta_X(t^a)}{\left[ w_\nu \right]} 
&=& 
\sum _{\nu' \in
  \mathcal A_{<a}(\nu)}  t^{\sum \nu' - \sum \nu}
K_{(<a)\nu'} 
\quad \text{(Prop.~\ref{P:Kbase}(a) and 
\eqref{eq:defK})} 
\end{eqnarray*}
 A little thought shows that if $\nu$ is the formalization of a partition,
and $\nu' \in \mathcal A_{<a}(\nu)$, then $m(\nu') \leq m(\nu)$, and 
if furthermore $m(\nu') = m(\nu)$, then $K_{(<a)\nu'} = K_{(<a)
  \nu}$.  Thus
\begin{eqnarray*}
\frac{\zzeta_X(t)}{\zzeta_X(t^a)}{[w_\nu]} &=&
\sum_{              \substack{\nu' \in  \mathcal A_{<a}(\nu) \\
m(\nu') < m(\nu)  }   }    
t^{\sum \nu' - \sum \nu}
K_{(<a)\nu'}  +
\sum _{\substack{\nu' \in
  \mathcal A_{<a}(\nu) \\ m(\nu') = m(\nu)}}  t^{\sum \nu' - \sum \nu}
K_{(<a)\nu'} 
\\
&=& \sum_{              \substack{\nu' \in  \mathcal A_{<a}(\nu) \\
m(\nu') < m(\nu)  }   }    
t^{\sum \nu' - \sum \nu}
K_{(<a)\nu'}  +
\left(\sum _{\substack{\nu' \in
  \mathcal A_{<a}(\nu) \\ m(\nu') = m(\nu)}}  t^{\sum \nu' - \sum
\nu} \right)
K_{(<a)\nu} 
\end{eqnarray*}
The result follows.
\qed

\epoint{Example (see \S \ref{promise59})}\label{E:distinctnu}
 If $\nu$ has all distinct elements greater than $1$, then
 Proposition~\ref{P:Krecursion} inductively yields
$$
K_{1^{\bullet} \nu}(t)= \frac {\zzeta_X(t)} {\zzeta_X\left(t^2\right)}
\frac{w_\nu }{\left(1+t\right)^{|\nu|}}.$$

Temporarily (for the purpose of Proposition~\ref{P:K}) define
$\mathcal S(\nu,a)$ to be the (finite) set of partitions $\mu$ with all parts at least $a$ such
that there exists a partition $\pi$ with elements each $<a$ with $
1^{|\nu|(a-1)} \nu \leq \pi\mu \not \geq 1^{|\nu|(a-1)-a}
av$.  In other words, these are the partitions (up to ``small parts''
$<a$) which can be obtained by merging in $|\nu| (a-1)$  ones with
$\nu$, but which cannot be obtained if $a$ of the ones are merged
together first.

\tpoint{Proposition} {\em \label{P:K} \quad
\newline  (a) 
%(degeneration formula for closed strata)
For a partition $\nu$ of positive integers, and an integer $a$ no bigger than
the smallest part of $\nu$,
$$
\left[ \wbar_{1^{j-a} a\nu} \right]= \left[ \wbar_{1^j \nu}  \right]-\sum_{\mu \in\mathcal
   S(\nu,a)} \left[ w_{(<a \vdash j -\sum \mu +\sum \nu  ) \mu} \right].
$$
\newline
(b)
%\tpoint{Proposition (recursion for $\Kbar_{1^{\bullet} \nu}(t)$)}\label{P:Kbarrecursion}
For any partition $\nu$ of positive integers all parts at least $a$,
 $$
\Kbar_{1^{\bullet} a\nu}(t)=\Kbar_{1^{\bullet} \nu}(t)t^{-a} - \sum_{\mu \in
  S(\nu,a)} K_{(<a)\mu}(t) t^{-a+\sum \mu-\sum \nu}.$$ 
  }

Theorem~\ref{T:Kform}(b) follows from 
inductively from Propositions~\ref{P:K}(b) and Theorem~\ref{T:Kform}(a),
and the base case
$\Kbar_{1^{\bullet} \varnothing} (t) =\zzeta_{X}(t)$.

\noindent {\em Proof.}
(a) 
By considering which $w_{\l}$  are contained in 
$\wbar_{1^j \nu}$ and in $\wbar_{1^{j-a} a\nu}$, 
we have
\begin{equation}
\left[ \wbar_{1^j \nu} \right]  -  \left[ \wbar_{1^{j-a} a\nu}
\right]=\sum_{1^j \nu \leq \lambda \not \geq 1^{j-a} a\nu} 
\left[ w_\lambda \right].\label{eq:T}
\end{equation}
We give a name to the partitions appearing on the right side of \eqref{eq:T}:
let \highlight{$\mathcal T (j):=\{\lambda\; | \; 1^j \nu \leq \lambda \not \geq 1^{j-a} a\nu \}$}.
For each $\lambda\in \mathcal T (j)$, we write  $\lam=b(\lam) s(\lam)$, where the ``big'' part
$b(\lam)$ is a partition composed of the elements of $\lam$ that are $\geq a$ and the ``small'' part
$s(\lam)$ is the rest.  Note that in any merge that created $\lam$
from $1^j \nu$, only $1$'s can contribute to the
$s(\lam)$ part.

Note that if all elements of an integer partition $\mu$ are at least
$a$, then (if $\sum 1^j \mu = \sum \lambda$)  whether $1^j\mu\leq \lambda$ depends only on $b(\lambda)$.
In particular, if $\mu=b(\lam)$ for some $\lam\in \mathcal T (j)$,
then for all partitions $\pi$ of $j-\sum \mu +\sum \nu$ into elements
less than $a$, we have that $\mu \pi\in \mathcal T (j)$.

Also note that
\begin{eqnarray}  \label{eq:ant} \big\{\mu \; | \;  \mu=b(\lam) \textrm{ for some }\lam\in  \mathcal T
(j) \big\}  &=& \big\{\mu \; | \; \mu  = b(\lam) \textrm{ for some }\lam\in \mathcal T
(|\nu|(a-1)) \\ & & \quad  \textrm{ and } \sum\mu- \sum \nu\leq j \big\}.\nonumber
\end{eqnarray}
%For all $\lambda\in \mathcal T (j)$, we have $\sum b(\lam)-\sum
%\nu\leq |\nu|(a-1)$.  Thus for $j\geq |\nu|(a-1)$, we can simplify
%\eqref{eq:ant} to
%$$\{\mu \; |  \; \mu=b(\lam) \textrm{ for some }\lam\in  \mathcal T
%(j)\}=\{\mu \; | \;  \mu=b(\lam) \textrm{ for some }\lam\in \mathcal T (|\nu|(a-1)) \}.$$

Thus
\begin{eqnarray*}
\left[ \wbar_{1^j \nu} \right]- \left[ \wbar_{1^{j-a} a\nu}  \right]&=& \sum_{\lambda\in \mathcal
  T (j)} \left[ w_\lambda \right]\quad  \text{(by definition of $\mathcal  T(j)$)}\\
&=&
\sum_{\mu \in S(\nu,a)} \left[ w_{(<a \vdash j-\sum\mu +\sum\nu )\mu} \right],
\end{eqnarray*}
where of course $[ w_{(<a \vdash j-\sum\mu +\sum\nu )\mu}]=0$ if $j-\sum\mu +\sum\nu<0.$

(b) Multiply both sides of (a)  by $t^{j-a}$, and sum over all $j$.  
\qed

\bpoint{Another example} Proposition~\ref{P:K}(b) gives a recursion to
compute $\Kbar_{1^{\bullet} \nu}$ in all cases, but for some $\nu$ we
have more efficient formulas.  One example is
Proposition~\ref{P:Kbase}(b) above.  Another is the following Lemma,
which provoked Conjecture~E (\S~\ref{c:q}), a topological conjecture about Betti
numbers.

\tpoint{Lemma} {\em \label{L:jksrecursion}
For $1<a\leq b$ and $r\geq 0$,
\begin{align*}
\wbar_{1^{j-a}ab^r} =  \wbar_{1^{j}b^r} - \wbar_{x^{j}y^r} +  \wbar_{x^{j-a}(ax)y^r},
\end{align*}
where $x$ and $y$ are formal variables. 
}

(Caution:  $a$ and $b$ are integers, while $x$ and $y$ are formal
variables --- this is key to the argument!) 

\noindent {\em Proof.}
If $\mu$ is a partition, let $\mathcal R_\mu$ be the set of partitions
$\geq \mu$, i.e. obtainable from $\mu$ by merging.
We have a map of posets $\mathcal R_{x^j y^r}\ra \mathcal R_{1^j b^r}$ sending $x\mapsto 1$ and $y\mapsto b$.  
%Note that $\mathcal R_{1^j b^r}\setminus \mathcal R_{1^{j-a}ab^r}$ consists of the partitions in $\mathcal R_{1^j b^r}$ for which when you take the standard residues of %the elements modulo $k$ they add (in $\Z$!) to $r$ (ignoring the formal elements from $v$).
We claim that the map $\mathcal R_{x^j y^r}\ra \mathcal R_{1^j b^r}$ restricts to a
{\em bijection} $\mathcal R_{x^j y^r}\setminus \mathcal R_{x^{j-a} (ax) y^r} \ra \mathcal R_{1^j b^r}\setminus
\mathcal R_{1^{j-a}ab^r}$, and that this bijection preserves the multiplicity sequence of each partition.  

First, we will see that $\mathcal R_{x^j y^r}\setminus \mathcal R_{x^{j-a} (ax) y^r}$ does map to $\mathcal R_{1^j b^r}\setminus \mathcal R_{1^{j-a}ab^r}$.
Consider an element $\mu$ of $\mathcal R_{x^j y^r}\setminus \mathcal R_{x^{j-a} (ax) y^r}$ that maps to $\lambda$ in $\mathcal R_{1^j b^r}$. 
 Each element of $\mu$ has at most $(a-1)$ $x$'s, and thus the reduction of an element of $\lambda$ modulo $b$
 is between $0$ and $a-1$.  In particular, the sum of these reductions (as integers, not modulo $b$)
 is $j$.  
If $\lambda$ were in $\mathcal R_{1^{j-a}ab^r}$, it would either have an element whose reduction is between $a$ and $b-1$ modulo $b$, 
 or the sum of the reductions modulo $b$ of the elements of $\lambda$
 would be less than $j$.

Second, given $\lambda \in \mathcal R_{1^j b^r}\setminus \mathcal
R_{1^{j-a}ab^r}$, looking at the residues of the elements modulo $b$,
we know where all the $1$'s have gone in any merge, and thus where all
the $b$'s are, determining a pre-image on $\mathcal R_{x^j y^r}
\setminus \mathcal R_{x^{j-a} (ax) y^r}$ uniquely.    Finally, since
no element of $\lambda\in \mathcal R_{x^j y^r}\setminus \mathcal
R_{x^{j-a} (ax) y^r}$ has more than $(a-1)$ $x's$, if two elements
$c_1x+c_2y=c_3x+c_4y$ ($c_1, \dots, c_4 \in \mathbb{Z}^{\geq 0}$) are equal after the map to $\mathcal R_{1^j b^r}$, then we have
$$
c_1+c_2b=c_3+c_4b
$$
for $0\leq c_1,c_3\leq a-1$, and thus $c_1=c_3$ and $c_2=c_4$.
This shows that  the map $\mathcal R_{x^j y^r}\setminus \mathcal R_{x^{j-a} (ax) y^r} \ra \mathcal R_{1^j b^r}\setminus \mathcal R_{1^{j-a}ab^r}$ preserves
multiplicity sequences.

The bijection $\mathcal R_{x^j y^r}\setminus \mathcal R_{x^{j-a} (ax) y^r} \ra \mathcal R_{1^j b^r}\setminus \mathcal R_{1^{j-a}ab^r}$ thus gives
$$\wbar_{1^jb^r}-\wbar_{1^{j-a}ab^r}=
\sum_{\lambda \in \mathcal R_{1^j b^r}\setminus \mathcal R_{1^{j-a}ab^r}} w_\lambda
 = \sum_{\mu \in \mathcal R_{x^j y^r}\setminus \mathcal R_{x^{j-a} (ax) y^r}} w_\mu =   \wbar_{x^{j}y^r} -  \wbar_{x^{j-a}(ax)y^r}.
$$
\qed

\tpoint{Proposition (see \S \ref{promise513} and Conjecture~\ref{c:q})} {\em \label{P:Kjksrecursion}
Given $1<a\leq b$ and $r\geq 0$, we have
$$
\Kbar_{1^{\bullet} ab^r} (t)
 =  \Kbar_{1^{\bullet} b^r} (t) t^{-a} - \frac {\zzeta_X(t)t^{-a}} {\zzeta_X(t^a)}
\left[  \Sym^r X \right].$$}

%\wbar_{1^{j-a}ab^r} =  \wbar_{1^{j}b^r} - \wbar_{x_1^{j}x_2^r} +  \wbar_{x_1^{j-a}(ax_1)x_2^r},

\noindent {\em Proof.}
Multiplying Lemma~\ref{L:jksrecursion} by $t^{j-a}$ and summing over $j$, we obtain
\begin{align*}
\Kbar_{1^{\bullet} ab^r} (t) &=  \Kbar_{1^{\bullet} b^r} (t) t^{-a} -\sum_j \left[ \Sym^j
  X \right] \left(  \left[ \Sym^r X \right] t^{j-a} \right) +\sum_j \wbar_{x_1^{j-a}(ax_1)}t^{j-a} \wbar_{x_2^r} \\
 &=  \Kbar_{1^{\bullet} b^r} (t) t^{-a} - \zzeta_X(t)   \left[ \Sym^r X \right] t^{-a}  
+\sum_j \left( \left[ \Sym^{j}X  \right] -  \left[ w_{(<a \vdash j)}
  \right] \right)t^{j-a}  \left[ \Sym^r X \right] \\
 &=  \Kbar_{1^{\bullet} b^r}(t) t^{-a} - \zzeta_X(t)   \left[ \Sym^r X \right] t^{-a}  +
\left( \zzeta_X(t) - K_{(<a)}(t) \right) \left[ \Sym^r X \right] t^{-a}.
\end{align*}
The result then follows from Proposition~\ref{P:Kbase}(a).
\qed

\epoint{Example (see \S \ref{promise513} and Conjecture~\ref{c:q})
}  
For $1<a\leq b$ and $r\geq 0$, we have
$$
\overline{K}_{1^{\bullet} a b^r}(t) = t^{-a-rb}\left( \zzeta_X(t)-
  \frac{\zzeta_X(t)}{\zzeta_X(t^{b})}\left(\sum_{i=0}^{r-1} {\left[
        \Sym^i X \right]}{t^{bi}}\right)
-\frac{\zzeta_X(t)}{\zzeta_X(t^{a})}
{ \left[ \Sym^r  X \right]}{t^{rb}}\right),
$$
by applying Proposition~\ref{P:Kjksrecursion} inductively.
(Note that $\overline{K}_{1^{\bullet} b^r}(t)$ should be 
interpreted as 
$\overline{K}_{1^{\bullet} b  b^{r-1}}(t)$ to be computed inductively.)
For example, if $X = \A^d$, then $\zzeta_X(t) =
1 / (1- \M t)$ (see the start of \S \ref{S:contract}) and
$\overline{K}_{1^{\bullet} a b^r}(t) =  M^{r+1} / (1-Mt)$.\label{ex:last}

%\bibliographystyle{amsplain}
%\bibliography{myrefs.bib}

\end{document}